
%

\input amstex.tex
\loadmsam
\loadmsbm
\loadbold
\input amssym.tex
\baselineskip=13pt plus 2pt
\documentstyle{amsppt}
\pageheight{45pc}
\pagewidth{33pc}

\magnification=1200
\overfullrule=0pt

\def\enddemos{\hfill$\square$\enddemo}
\def\hb{\hfil\break}
\def\n{\noindent}
\def\ov{\overline}
\def\smatrix{\smallmatrix}

\def\pmatrix{\left(\smatrix}
\def\endpmatrix{\endsmallmatrix\right)}
\def\upi{\pmb{\pi}}
\def\pii{\pmb{\pi}}
\def\upsi{\pmb{\psi}}

\def\r{\bold r}
\def\s{\bold s}

\def\A{\Bbb A}
\def\C{\Bbb C}
\def\F{\Bbb F}
\def\G{\Bbb G}

\def\R{\Bbb R}
\def\Z{\Bbb Z}

\def\1{\bold 1}
\def\v{\bold v}
\def\w{\bold w}
\def\uC{\bold C}
\def\uG{\bold G}
\def\uT{\bold T}
\def\Gal{\operatorname{Gal}}
\def\max{\operatorname{max}}
\def\Vol{\operatorname{vol}}

\def\GL{\operatorname{GL}}
\def\Gp{\operatorname{Gp}}
\def\SL{\operatorname{SL}}
\def\PGL{\operatorname{PGL}}
\def\PGp{\operatorname{PGp}}
\def\SO{\operatorname{SO}}

\def\Ind{\operatorname{Ind}}
\def\tr{\operatorname{tr}}
\def\Re{\operatorname{Re}}
\def\Int{\operatorname{Int}}
\def\Lie{\operatorname{Lie}}
\def\Ad{\operatorname{Ad}}

\def\diag{\operatorname{diag}}

\leftheadtext{Yuval Z. Flicker and Dmitrii Zinoviev}
\rightheadtext{Twisted character of a small representation}
\topmatter
\title Twisted character of a small representation of PGL(4)\endtitle
\author Yuval Z. Flicker and Dmitrii Zinoviev\endauthor
\footnote"~"{\n Department of Mathematics, The Ohio State University,
231 W. 18th Ave., Columbus, OH 43210-1174;\hb
Institute of Mathematics, The Hebrew University, Givat Ram, Jerusalem 91904, 
Israel;\hb email: flicker\@math.ohio-state.edu.\hb
\indent Institute for Information Transmission Problems, Russian Academy 
of Sciences, Moscow, Russia;\hb
email: zinov\@iitp.ru.\hb 
\n Keywords: Representations of $p$-adic groups, explicit character
computations, twisted endoscopy, transpose-inverse twisting, instability,
liftings\hb
\indent 2000 Mathematics Subject Classification: 22E55, 11F70, 11F85, 11F46,
20G25, 22E35.}

\abstract We compute by a purely local method the elliptic $\theta$-twisted 
character $\chi_\pi$ of the representation $\pi=I_{(3,1)}(1_3)$ of 
$\PGL(4,F)$. Here $F$ is a $p$-adic field; $\theta$ is the 
``transpose-inverse'' automorphism of $G=\PGL(4,F)$; $\pi$ is the 
representation of $\PGL(4,F)$ normalizedly induced from the trivial 
representation of the maximal parabolic subgroup of type $(3,1)$. 
Put $\uC=\{(g_1,g_2)\in \GL(2)\times\GL(2);\,\det(g_1)=\det(g_2)\}/\G_m$
($\G_m$ embeds diagonally). It is a $\theta$-twisted elliptic endoscopic 
group of $\PGL(4)$. We deduce from the computation that $\chi_\pi$ is an 
unstable function: its value at one twisted regular elliptic conjugacy 
class with norm in $C=\uC(F)$ is minus its value at the other class 
within the twisted stable conjugacy class, and $0$ at the classes 
without norm in $C$. Moreover $\pi$ is the unstable endoscopic lift 
of the trivial representation of $C$.

Naturally, this computation plays a role in the theory of lifting
from $\uC$(=``SO(4)'') and $\PGp(2)$ to $\uG=\PGL(4)$ using the 
trace formula, to be discussed elsewhere (see [F$'$]).

Our work develops a 4-dimensional analogue of the model of the small 
representation of $\PGL(3,F)$ introduced with Kazhdan in [FK] in a
3-dimensional case. It uses the classification of twisted stable and
unstable regular conjugacy classes in $\PGL(4,F)$ of [F], motivated by 
Weissauer [W]. It extends the local method of computation introduced by 
us in [FZ]. An extension of our work here to apply to similar 
representations of GL(4,$F$) whose central character is non trivial
has recently been given in [FZ$'$].

\endabstract
\endtopmatter

\document
\heading Introduction\endheading
Let $\pi$ be an admissible representation (see Bernstein-Zelevinsky [BZ], 
2.1) of a $p$-adic reductive group $G$. Its character $\chi_\pi$ is a 
complex valued function defined by $\tr\pi(fdg)=\int_G\chi_\pi(g)f(g)dg$ 
for all complex valued smooth compactly supported measures $fdg$ ([BZ],
2.17). It is smooth on the regular set of the group $G$. The character
is important since it characterizes the representation up to equivalence.
A fundamental result of Harish-Chandra [H] establishes that the character 
is a locally integrable function in characteristic zero. 

Let $\theta$ be an automorphism of finite order of the group $G$. Define
${}^\theta\pi$ by ${}^\theta\pi(g)=\pi(\theta(g))$. When $\pi$ is invariant 
under the action of $\theta$ (thus ${}^\theta\pi$ is equivalent to $\pi$), 
Shintani and others introduced an extension of $\pi$ to the semidirect 
product $G\rtimes\langle\theta\rangle$. The twisted character 
$\chi_\pi(g\times\theta)$ is defined by $\tr\pi(fdg\times\theta)
=\int_G\chi_\pi(g\times\theta)f(g)dg$ for all $fdg$. It depends only 
on the $\theta$-conjugacy class $\{hg\theta(h)^{-1};h\in G\}$ of $g$. 
It is again smooth on the $\theta$-regular set, and characterizes the 
$\theta$-invariant irreducible $\pi$ up to isomorphism. Moreover, it
is locally integrable (see Clozel [C]) in characteristic zero.

Characters provide a very precise tool to express a relation of 
representations of different groups, called lifting. It was studied 
extensively by Shintani and others in the case of base change. It was 
studied also in non base change situations such as twisting by characters
(Kazhdan [K], Waldspurger [Wa]), and the symmetric square lifting
from SL(2) to PGL(3) ([Fsym], [FK]). In this last case twisted 
characters of $\theta$-invariant representations of PGL(3) are
related to packets of representations of SL(2), and $\theta$ is
the involution sending $g$ to its transpose-inverse.

{\it The aim of the present work is to compute the twisted $($by $\theta)$
character of a specific representation $\pi=I_{(3,1)}(1_3)$, of the
group $G=\PGL(4,F)$, $F$ a $p$-adic field, $p$ odd}. This $\pi$ is 
normalizedly induced from the trivial representation $1_3$ of the 
standard (upper triangular) maximal parabolic subgroup $P$ of type 
$(3,1)$. It is invariant under the involution 
$\theta(g)=J^{-1}{}^tg^{-1}J$, where
$J=\pmatrix 0 &w\\-w & 0\endpmatrix$ and $w=\pmatrix 0 &1\\1 & 0\endpmatrix$.

A natural setting for the statement of our result is the theory of
liftings to the group $\uG=\PGL(4)$ from its $\theta$-twisted
endoscopic group (see Kottwitz-Shelstad [KS])
$$\uC=\{(g,g')\in\GL(2)\times\GL(2);\,\det g=\det g'\}/\G_m.$$
Here the multiplicative group $\G_m=\GL(1)$
embeds as $z\mapsto(zI_2,zI_2)$, $I_2$ is the identity $2\times 2$
matrix. The corresponding map $\lambda_1$ of dual groups is simply
the natural embedding in $\hat{G}=\SL(4,\C)$ of
$\hat{C}=Z_{\hat{G}}(\hat{s}\hat\theta)=``\SO(4,\C)$''
$$=\left\{g\in\hat{G}=\SL(4,\C);\,g\hat{s}J{}^tg=\hat{s}J=
\pmatrix 0 &\omega\\ \omega^{-1} & 0\endpmatrix\right\}
=\SO\left(\pmatrix 0 &\omega\\ \omega^{-1} & 0\endpmatrix,\C\right)$$
$$=\left\{\pmatrix aB & bB\\ cB & dB\endpmatrix;\,
\left(A=\pmatrix a & b\\ c & d\endpmatrix,B\right)
\in(\GL(2,\C)\times\GL(2,\C),\,\,\det A\cdot\det B=1)/\C^\times\right\}.$$
Here $z\in\C^\times$ embeds as the central element $(z,z^{-1})$, and
$\hat{s}=\diag(-1,1,-1,1)$ and $\omega=\pmatrix 0&-1\\ 1&0\endpmatrix$.
Thus $\hat{C}$ is the $\hat\theta$-centralizer in $\hat{G}$ of the
semisimple element $\hat{s}$, and $\hat\theta$ is defined on $\hat{G}$
by the same formula that defines $\theta$ on $G$.

Our result can be viewed as asserting that the $\theta$-invariant
representation $\pi$ of $G=\uG(F)$ is the endoscopic lift of the trivial
representation of $C=\uC(F)$. To state this we note that the embedding
$\lambda_1:\hat{C}\to\hat{G}$ defines a norm map. This norm map relates
the stable $\theta$-conjugacy classes in $G$ with stable conjugacy classes
in $C$. A stable conjugacy class is the intersection $G\cap(\Int(\uG(\ov{F}))
(\gamma))$, for some $\gamma\in G$. The crucial case of the character
computation is that of $\theta$-elliptic elements. A stable 
$\theta$-conjugacy class consists of several $\theta$-conjugacy classes. 
The stable $\theta$-conjugacy classes of elements in $G$, and the 
$\theta$-conjugacy classes within the stable $\theta$-classes, have been 
described recently in [F], in analogy with the description of the conjugacy 
classes and the stable classes in the group of symplectic similitudes 
$\Gp(2,F)$ of Weissauer [W]. In fact in [F] we deal with $\theta$-classes 
in $\GL(4,F)$ while here we deal with the simpler case of $\PGL(4,F)$, so 
we give here full details of the description in our case.

There are four types of $\theta$-elliptic elements of $G$, named in [F], 
p. 16, and here (see the next section) I, II, III, IV, depending on their 
splitting behaviour. As in [F], our work relies on an explicit presentation 
of representatives of the $\theta$-conjugacy classes within the stable such 
classes in $G$, except that here we present a better looking set of such 
representatives.

The norm map, which we describe explicitly here, relates $\theta$-conjugacy
classes of types I and III to conjugacy classes in $C$. It does not relate
classes of types II, IV to classes in $C$.

{\it We prove that the $\theta$-character of $\pi$, $\chi_\pi(g\times\theta)$,
vanishes on $\theta$-regular elements $g$ of type II and IV}. The stable
$\theta$-conjugacy classes of types I and III come associated with a
quadratic extension $E/F$ in type I and $E/E_3$ in type III (in this case
$E_3/F$ is a quadratic extension, and $E/F$ is biquadratic). The two
$\theta$-conjugacy classes $g_r$ within the stable $\theta$-classes are
parametrized by $r$ in $F^\times/N_{E/F}E^\times$ in type I and by 
$E_3^\times/N_{E/E_3}E^\times$ in type III. {\it We show that the value 
of $\chi_\pi(g_r\times\theta)$, multiplied by a suitable Jacobian
$\Delta(g_r\theta)/\Delta_C(Ng)$, is $2\kappa(r)$}. Here $\kappa$
is the nontrivial character of $F^\times/N_{E/F}E^\times$ in type I 
and of $E_3^\times/N_{E/E_3}E^\times$ in type III.

In particular the character $\chi_\pi(g\times\theta)$ is an unstable
function, namely its value at one $\theta$-conjugacy class within a
stable $\theta$-conjugacy class of type I or III is negative its value
at the other $\theta$-conjugacy class.

Our result is a special case of the lifting with respect to $\lambda_1$
to the group $G=\PGL(4,F)$ of representations of the group $C=(\GL(2,F)
\times\GL(2,F))'/F^\times$, where the prime indicates $\det g=\det g'$
for the two components $(g,g')$, and $F^\times$ embeds diagonally. This
lifting is established in [F$'$] by means of a comparison of trace formulae,
the fundamental lemma of [F], and character relations, for generic and
nongeneric representations. It can be viewed as associating to a pair
$\pi_1$, $\pi_2$ of representations of $\GL(2,F)$ (the product of whose
central characters is 1) a product representation $\pi=\pi_1\boxtimes\pi_2$,
or $\lambda_1(\pi_1\times\pi_2)$, of $\PGL(4,F)$.

The case that we consider here is that where $\pi_1$ and $\pi_2$ are the
trivial representations of $\GL(2,F)$. Their product via the lifting 
$\lambda_1$, $1\boxtimes 1$ or $\lambda_1(1\times 1)$, is our 
$\pi=I_{(3,1)}(1_3)$.

This notion of multiplication should not be confused with that of induction
from the standard parabolic subgroup of type $(2,2)$ and $\pi_1\otimes\pi_2$ 
on its Levi factor, thus: $I_{(2,2)}(\pi_1\otimes\pi_2)$, which plays the
role of addition of $\pi_1$ and $\pi_2$. This multiplication is interesting
in particular since conjecturally the notions of multiplication 
$\pi_1\boxtimes\pi_2=\lambda_1(\pi_1\times\pi_2)$ and addition
$\pi_1\boxplus\pi_2=I_{(2,2)}(\pi_1\otimes\pi_2)$ -- suitably extended to
all GL($n$) -- give a Tannakian structure on the category of algebraic (=
smooth) representations of all the $\GL(n,F)$, and its motivic Galois group 
(see Deligne-Milne [DM]) plays a key role in the principle of functoriality.

However, the proof of [F$'$], based on trace formulae comparison and the
fundamental lemma of [F], is very involved. The current paper grew from an
attempt to provide a purely local and self contained proof of a key initial 
case, where the representations involved are not generic, where all principal
features can be explicitly viewed: the analysis of the $\theta$-conjugacy
classes, the norm map, the Jacobian factors and the transfer factors, and
the character relations can be computed directly to verify them without
relying on long and complex theories. 

This gives an independent verification of results obtainable by global 
techniques, by purely local and essentially elementary techniques.

Our method here is based on using a novel model of our representation
$\pi=I_{(3,1)}(1_3)$, different from the standard model of a parabolically
induced representation. It is a four dimensional analogue of a three
dimensional model introduced and used with Kazhdan in [FK] to compute 
the twisted by transpose-inverse character of the representation 
$\pi_3=I_{(2,1)}(1_2)$ of $\PGL(3,F)$ normalizedly induced from the trivial 
representation of the maximal parabolic subgroup. The original interest
of [FK] was in a case of the fundamental lemma for the symmetric square 
lifting. But a purely local and simpler proof was given later in [Fsym; Unit
elements]. In our case the fundamental lemma is established in [F]. 

The work of [FK] uses local arguments to 
compute the twisted character of $\pi_3$ on one of the two twisted 
conjugacy classes within the stable one (where the quadratic form is
anisotropic), and global arguments to reduce the computation on the other 
class (where the quadratic form is isotropic) to that computed by local 
means. A purely local computation for the second class is given in [FZ]. 
Here we develop this local computation in our four dimensional case. A
global type of argument as in [FK] is harder to apply as there are not
enough anisotropic quadratic forms in our case. Anyway, here we give a 
simpler, local proof.

We believe that our method of computation is applicable in many cases of
character computations, giving rise to a new theory of integration of
functions on $p$-adic domains, and we plan to return to this topic in
future work. See [FZ$'$] for an analogue of the present work in the 
case of GL(4,$F$) when the central character is non trivial.

We work only with a $p$-adic field $F$. However, the model of our 
representation makes sense also in the case of a real base field $F=\R$. 
We propose the question of verifying our formula for the twisted character 
of our $\pi$ also for the field $\R$. Of course in this case there is only 
one type of stable elliptic $\theta$-conjugacy class, namely type I, as the 
only algebraic field extension of the reals $\R$ is the field $\C$ of 
complex numbers. See the remark before Theorem I below.

We are deeply grateful to the referee for careful reading of this work.

\heading Conjugacy classes\endheading

Let $F$ be a local nonarchimedean field, and $R$ its ring of integers.
Put $\uG=\PGL(4)$, $G=\uG(F)$ and $K=\uG(R)$. Put 
$\uC=\{(g_1,g_2)\in \GL(2)\times \GL(2);\,\det(g_1)=\det(g_2)\}/\G_m$
($\G_m$ embeds diagonally),
$C=\uC(F)=\{(g_1,g_2)\in \GL(2,F)\times\GL(2,F);\,
\det(g_1)=\det(g_2)\}/F^\times$ and $K_C=\uC(R)$.
Put $J=(a_i\delta_{i,5-j})$, $a_1=a_2=1$, $a_3=a_4=-1$, and set 
$\theta(\delta)=J^{-1}{}^t\delta^{-1}J$ for $\delta$ in $G$.
Fix a separable algebraic closure $\ov{F}$ of $F$. 
The elements $\delta$, $\delta'$ of $G$ are called ({\it stably}) 
$\theta$-{\it conjugate} if there is $g$ in $G$ (resp. $\PGL(4,\ov{F}))$
with $\delta'=g^{-1}\delta\theta(g)$.

We recall some results of [F] concerning (stable) $\theta$-twisted regular
conjugacy classes. There are four types of $\theta$-elliptic classes,
but the norm map $N$ from $G$ to $C$ relates only the twisted classes
in $G$ of type I and III to conjugacy classes in $C$. We should then 
expect the twisted character of the representation considered here to 
vanish on the twisted classes of type II and IV.

A set of representatives for the $\theta$-conjugacy classes 
within a stable semisimple $\theta$-conjugacy class of type 
I in $\GL(4,F)$ which splits over a quadratic extension 
$E=F(\sqrt D)$ of $F$, $D\in F-F^2$, is parametrized by 
$(\r,\s)\in F^\times/N_{E/F}E^\times\times F^\times/N_{E/F}E^\times$ 
([F], p. 16). Representatives for the $\theta$-regular (thus $t\theta(t)$
is regular) stable $\theta$-conjugacy classes of type I in $\GL(4,F)$
which split over $E$ can be found in a torus $T=\uT(F)$, 
$\uT=h^{-1}\uT^\ast h$, $\uT^\ast$ denoting the diagonal subgroup 
in $\uG$, $h=\theta(h)$, and
$$T=\left\{t=\pmatrix a_1 & 0 & 0 & a_2D\\ 0 & b_1 & b_2D & 0 
\\ 0 & b_2 & b_1& 0\\ a_2 & 0 & 0 & a_1\endpmatrix
=h^{-1}t^\ast h;\quad t^\ast=\diag(a,b,\sigma b,\sigma a)\in T^\ast
\right\}.$$
Here $a=a_1+a_2\sqrt D,\,b=b_1+b_2\sqrt D\in E^\times$, and $t$ is
regular if $a/\sigma a$ and $b/\sigma b$ are distinct and not equal to
$\pm 1$. Note that here $T^\ast=\uT^\ast(F)$ where the Galois action
is that obtained from the Galois action on $T$.

A set of representatives for the $\theta$-conjugacy classes within
a stable $\theta$-conjugacy class can be chosen in $T$. Indeed, if
$t=h^{-1}t^\ast h$ and $t_1=h^{-1}t_1^\ast h$ in $T$ are stably
$\theta$-conjugate, then there is $g=h^{-1}\mu h$ with 
$t_1=gt\theta(g)^{-1}$, thus $t_1^\ast=\mu t^\ast\theta(\mu)^{-1}$
and $t_1^\ast\theta(t_1^\ast)=\mu t^\ast\theta(t^\ast)\mu^{-1}$.
Since $t$ is $\theta$-regular, $\mu$ lies in the $\theta$-normalizer of 
$\uT^\ast(\ov{F})$ in $\uG(\ov{F})$. Since the group $W^\theta(\uT^\ast,\uG)$
$=N^\theta(\uT^\ast,\uG)/\uT^\ast$, quotient by $\uT^\ast(\ov{F})$ of the
$\theta$-normalizer of $\uT^\ast(\ov{F})$ in $\uG(\ov{F})$, is represented
by the group $W^\theta(T^\ast,G)=N^\theta(T^\ast,G)/T^\ast$, quotient by 
$T^\ast$ of the $\theta$-normalizer of $T^\ast$ in $G$, we may modify $\mu$
by an element of $W^\theta(T^\ast,G)$, that is replace $t_1$ by a 
$\theta$-conjugate element, and assume that $\mu$ lies in $\uT^\ast(\ov{F})$.
In this case $\mu\theta(\mu)^{-1}=\diag(u,u',\sigma u',\sigma u)$ (since $t$, 
$t_1$ lie in $T^\ast$), with $u=\sigma u$, $u'=\sigma u'$ in $F^\times$. 
Such $t$, $t_1$ are $\theta$-conjugate if $g\in G$, thus $g\in T$, so 
$\mu=\diag(v,v',\sigma v',\sigma v)\in T^\ast$ and $\mu\theta(\mu)^{-1}
=\diag(v\sigma v,v'\sigma v',v'\sigma v',v\sigma v)$. Hence a set
of representatives for the $\theta$-conjugacy classes within the
stable $\theta$-conjugacy class of the $\theta$-regular $t$ in $T$
is given by $t\cdot\diag(\r,\s,\s,\r)$, where 
$\r,\,\s\in F^\times/N_{E/F}E^\times$.
Clearly in $\PGL(4,F)$ the $\theta$-classes within a stable class 
are parametrized only by $\r$, or equivalently only by $\s$.

A set of representatives for the $\theta$-conjugacy classes 
within a stable semisimple $\theta$-conjugacy class of type 
II in $\GL(4,F)$ which splits over the biquadratic extension
$E=E_1E_2$ of $F$ with Galois group $\langle\sigma,\tau\rangle$, 
where $E_1=F(\sqrt D)=E^\tau$, $E_2=F(\sqrt{AD})=E^{\sigma\tau}$,
$E_3=F(\sqrt A)=E^\sigma$ are quadratic extensions of $F$, thus 
$A,D\in F-F^2$, is parametrized by $\r\in F^\times/N_{E_1/F}E_1^\times$, 
$\s\in F^\times/N_{E_2/F}E_2^\times$ ([F], p. 16). It is given by
$$\pmatrix a_1\r & 0 & 0 & a_2D\r\\ 0 & b_1\s & b_2AD\s & 0\\ 
0 & b_2\s & b_1\s& 0\\ a_2\r & 0 & 0 & a_1\r\endpmatrix
=h^{-1}t^\ast h\cdot\diag(\r,\s,\s,\r),\qquad 
t^\ast=\diag(a,b,\tau b,\sigma a).$$
Here $a=a_1+a_2\sqrt D\in E_1^\times$, $b=b_1+b_2\sqrt{AD}\in E_2^\times$, 
$\theta(h)=h$. In $\PGL(4,F)$ the $\theta$-classes within a stable class are 
parametrized only by $\r$, or equivalently only by $\s$.

A set of representatives for the $\theta$-conjugacy classes 
within a stable semisimple $\theta$-conjugacy class of type 
III in $\GL(4,F)$ which splits over the biquadratic extension
$E=E_1E_2$ of $F$ with Galois group $\langle\sigma,\tau\rangle$, 
where $E_1=F(\sqrt D)=E^\tau$, $E_2=F(\sqrt{AD})=E^{\sigma\tau}$,
$E_3=F(\sqrt A)=E^\sigma$ are quadratic extensions of $F$, thus 
$A,D\in F-F^2$, is parametrized by 
$r(=r_1+r_2\sqrt A)\in E_3^\times/N_{E/{E_3}}E^\times$ 
([F], p. 16). Representatives for the stable regular $\theta$-conjugacy
classes can be taken in the torus $T=h^{-1}T^\ast h$, consisting of
$$t=\pmatrix {\pmb a} & {\pmb b}D\\ 
{\pmb b} & {\pmb a}\endpmatrix=h^{-1}t^\ast h,\qquad 
t^\ast=\diag(\alpha,\tau\alpha,\sigma\tau\alpha,\sigma\alpha),$$
where $h=\theta(h)$ is described in [F], p. 16. This $t$ is 
$\theta$-regular when $\alpha/\sigma\alpha$, $\tau(\alpha/\sigma\alpha)$
are distinct and $\not=\pm 1$. Here
$$\qquad{\pmb a}=\pmatrix a_1 & a_2A\\ a_2 & a_1\endpmatrix,
\qquad{\pmb b}=\pmatrix b_1 & b_2A\\ b_2 & b_1\endpmatrix;
\qquad{\text{put also}}
\qquad{\pmb r}=\pmatrix r_1 & r_2A\\ r_2 & r_1\endpmatrix.$$
Further $\alpha=a+b\sqrt D\in E^\times$, $a=a_1+a_2\sqrt A\in E_3^\times$,
$b=b_1+b_2\sqrt A\in E_3^\times$, $\sigma\alpha=a-b\sqrt D$,
$\tau\alpha=\tau a+\tau b\sqrt D$. Representatives for all 
$\theta$-conjugacy classes within the stable $\theta$-conjugacy class 
of $t$ can be taken in $T$. In fact if $t'=gt\theta(g)^{-1}$ lies in $T$
and $g=h^{-1}\mu h$, $\mu\in\uT^\ast(\ov{F})$, then $\mu\theta(\mu)^{-1}
=\diag(u,\tau u,\sigma\tau u,\sigma u)$ has $u=\sigma u$,
thus $u\in E_3^\times$. If $g\in T$, thus $\mu\in T^\ast$, then
$\mu=\diag(v,\tau v,\sigma\tau v,\sigma v)$ and
$\mu\theta(\mu)^{-1}=\diag(v\sigma v,\tau v\sigma\tau v,
\tau v\sigma\tau v,v\sigma v)$, with 
$v\sigma v\in N_{E/E_3}E^\times$. We conclude that a set of 
representatives for the $\theta$-conjugacy classes within the stable
$\theta$-conjugacy class of $t$ is given by $t\cdot\diag(\pmb r,\pmb r)$,
$r\in E_3^\times/N_{E/E_3}E^\times$.

Representatives for the stable regular $\theta$-conjugacy classes of type 
IV can be taken in the torus $T=h^{-1}T^\ast h$, consisting of
$$t=\pmatrix {\pmb a} & {\pmb b}{\pmb D}\\ 
{\pmb b} & {\pmb a}\endpmatrix=h^{-1}t^\ast h,\qquad 
t^\ast=\diag(\alpha,\sigma\alpha,\sigma^3\alpha,\sigma^2\alpha),$$
where $h=\theta(h)$ is described in [F], p. 18. Here $\alpha$ ranges
over a quadratic extension $E=F(\sqrt D)=E_3(\sqrt D)$ of a quadratic
extension $E_3=F(\sqrt A)$ of $F$. Thus $A\in F-F^2$, $D=d_1+d_2\sqrt A$
lies in $E_3-E_3^2$ where $d_i\in F$. The normal closure $E'$ of $E$
over $F$ is $E$ if $E/F$ is cyclic with Galois group $\Z/4$, or a
quadratic extension of $E$, generated by a fourth root of unity $\zeta$,
in which case the Galois group is the dihedral group $D_4$. In both cases
the Galois group contains an element $\sigma$ with $\sigma\sqrt A=
-\sqrt A$, $\sigma\sqrt D=\sqrt{\sigma D}$, $\sigma^2\sqrt D=-\sqrt D$.
In the $D_4$ case $\Gal(E'/F)$ contains also $\tau$ with $\tau\zeta=-\zeta$,
we may choose $D=\sqrt A$, $\tau D=D$ and $\sigma\sqrt D=\zeta\sqrt D$.

In any case, $t$ is $\theta$-regular when $\alpha\not=\sigma^2\alpha$.
We write $\alpha=a+b\sqrt D\in E^\times$, $a=a_1+a_2\sqrt A\in E_3^\times$,
$b=b_1+b_2\sqrt A\in E_3^\times$, 
$\sigma\alpha=\sigma a+\sigma b\sqrt{\sigma D}$,
$\sigma^2\alpha=a-b\sqrt D$. Also
$$\qquad{\pmb a}=\pmatrix a_1 & a_2A\\ a_2 & a_1\endpmatrix,
\qquad{\pmb b}=\pmatrix b_1 & b_2A\\ b_2 & b_1\endpmatrix,
\qquad{\pmb D}=\pmatrix d_1 & d_2A\\ d_2 & d_1\endpmatrix.$$
Representatives for all $\theta$-conjugacy classes within the stable 
$\theta$-conjugacy class of $t$ can be taken in $T$. In fact if 
$t'=gt\theta(g)^{-1}$ lies in $T$ and $g=h^{-1}\mu h$, 
$\mu\in\uT^\ast(\ov{F})$, then $\mu\theta(\mu)^{-1}
=\diag(u,\sigma u,\sigma^3 u,\sigma^2 u)$ has $u=\sigma^2 u$,
thus $u\in E_3^\times$. If $g\in T$, thus $\mu\in T^\ast$, then
$\mu=\diag(v,\sigma v,\sigma^3 v,\sigma^2 v)$ and
$\mu\theta(\mu)^{-1}=\diag(v\sigma^2 v,\sigma(v\sigma^2 v),
\sigma(v\sigma^2 v),v\sigma^2 v)$, with 
$v\sigma v\in N_{E/E_3}E^\times$. It follows that a set of 
representatives for the $\theta$-conjugacy classes within the stable
$\theta$-conjugacy class of $t=h^{-1}t^\ast h
=\pmatrix {\pmb a} & {\pmb b}{\pmb D}\\ 
{\pmb b} & {\pmb a}\endpmatrix$, where
$t^\ast=\diag(\alpha,\sigma\alpha,\sigma^3\alpha,\sigma^2\alpha)$, is
given by multiplying $\alpha$ by $r$, that is $t^\ast$ by 
$t_0^\ast=\diag(r,\sigma r,\sigma^3 r,\sigma^2 r)$, where $r=\sigma^2 r$ 
ranges over a set of representatives for $E_3^\times/N_{E/E_3}E^\times$.
Now $t_0=h^{-1}t_0^\ast h=\pmatrix {\pmb r}&\\&{\pmb r}\endpmatrix$.
Hence a set of representatives is given by $t\cdot\diag(\pmb r,\pmb r)$,
$r\in E_3^\times/N_{E/E_3}E^\times$.

\heading Norm map\endheading

The norm map $N:\uG\to\uC$ is defined on the diagonal torus $\uT^\ast$
of $\uG$ by $$N(\diag(a,b,c,d))=(\diag(ab,cd),\diag(ac,bd)).$$ Since
both components have determinant $abcd$, the image of $N$
is indeed in $\uC$. 

In type I we have $a,b\in E^\times=F(\sqrt D)^\times$ and the norm map becomes
$$N(\diag(a,b,\sigma b,\sigma a))=(\diag(ab,\sigma(ab)),\diag(a\sigma b,
b\sigma a)).$$

In type III we have $\alpha\in E^\times$, $\alpha\tau\alpha\in E_1^\times$,
$\alpha\sigma\tau\alpha\in E_2^\times$, and the norm map becomes
$$N(\diag(\alpha,\tau\alpha,\sigma\tau\alpha,\sigma\alpha))=
(\diag(\alpha\tau\alpha,\sigma\alpha\sigma\tau\alpha),
\diag(\alpha\sigma\tau\alpha,\tau\alpha\sigma\alpha)).$$

The $\diag(\ast,\ast)$ define conjugacy classes in $\GL(2,F)$, and 
since both components of $N(\ast)$ have equal determinants in $F^\times$, 
the norm map defines a conjugacy class in $C=\uC(F)$ for each $\theta$-stable
conjugacy class of type I or III in $G=\uG(F)$.

In types II and IV no conjugacy class in $C$ corresponds
to the image of the map $N$.

\heading Jacobians\endheading

The character relation that we study relates the product of the value
at $t$ of the twisted character of our representation $\pi=I_{(3,1)}(\pi_1)$ 
by a factor $\Delta(t\times\theta)$, with the product by a factor 
$\Delta_C(Nt)$ of the value at $Nt$ of the character of the representation 
$\pi_C$ of $C$ which lifts to $\pi$.

The factor $\Delta(t\times\theta)$ is defined by
$$\Delta(t\times\theta)^2=|\det(1-\Ad(t\theta))|\Lie(G/T)|.$$
Here $t$ lies in the $\theta$-invariant torus $T$ which we take
to have the form $T=h^{-1}T^\ast h$, $T^\ast$ is the diagonal
subgroup and $h=\theta(h)$. Thus in the formula for
$\Delta(t\times\theta)$ we may replace $t=h^{-1}t^\ast h$ and $T$ 
by the diagonal $t^\ast$ and $T^\ast$. Note that $\Lie(G/T^\ast)
=\Lie U\oplus\Lie U^-$, and the upper and lower triangular subgroups
$U$, $U^-$ are $\theta$-invariant. We have
$$|\det(1-\Ad(t\theta))|\Lie U|=|\prod_\Theta 
(1-\sum_{\beta\in\Theta}\beta(t))|,$$
where the product ranges over the orbits $\Theta$ of $\theta$ in 
the set of positive roots $\beta>0$, and the sum ranges over the roots 
in the $\theta$-orbit. Thus for $t=\diag(a,b,c,d)$ we obtain
$$\left|\left(1-{a\over b}{c\over d}\right)\left(1-{a\over c}{b\over d}\right)
\left(1-{a\over d}\right)\left(1-{b\over c}\right)\right|.$$
Further,
$$|\det(1-\Ad(t\theta))|\Lie U^-|=\delta(t\theta)^{-1}
|\det(1-\Ad(t\theta))|\Lie U|$$
where
$$\delta(t\theta)=\prod_\Theta((\sum_{\beta\in\Theta}\beta)(t))
=({a\over b}{c\over d})({a\over c}{b\over d})({b\over c})({a\over d})
=\prod_{\beta>0}\beta(t)=\delta(t).$$
Altogether
$$\Delta(t\theta)=\left|{{(ac-bd)^2}\over{abcd}}\cdot{{(ab-cd)^2}\over{abcd}}\cdot
{{(a-d)^2}\over{ad}}{{(b-c)^2}\over{bc}}\right|^{1/2}.$$

Similarly,
$$\Delta_C(Nt)=\delta_C^{-1/2}(Nt)|\det(1-Nt)|\Lie U_C|=
\left|{ab\over{cd}}\cdot{ac\over{bd}}\right|^{-1/2}
\left|\left(1-{ab\over{cd}}\right)\left(1-{ac\over{bd}}\right)\right|,$$
and so 
$${\Delta(t\theta)\over\Delta_C(Nt)}=\left|{(a-d)^2\over{ad}}\cdot
{(b-c)^2\over{bc}}\right|^{1/2}.$$

Then in case I if $t=\diag(a,b,\sigma b,\sigma a)$, $a=a_1+a_2\sqrt D$,
$b=b_1+b_2\sqrt D$, we get 
$${\Delta(t\theta)\over\Delta_C(Nt)}
=\left|{(a-\sigma a)^2\over{a\sigma a}}\cdot
{(b-\sigma b)^2\over{b\sigma b}}\right|^{1/2}
=\left|{(2a_2\sqrt D)^2\over{a_1^2-a_2^2D}}
\cdot{(2b_2\sqrt D)^2\over{b_1^2-b_2^2D}}\right|^{1/2}.$$

In case III, if $t=\diag(\alpha,\tau\alpha,\sigma\tau\alpha,\sigma\alpha)$,
$\alpha=a+b\sqrt D$, $a=a_1+a_2\sqrt A$, $b=b_1+b_2\sqrt A$, 
$\sigma\alpha=a-b\sqrt D$, $\tau\alpha=\tau a+\tau b\sqrt D$, 
$\alpha-\sigma\alpha=2b\sqrt D$, $\tau(\alpha-\sigma\alpha)=2\tau b\sqrt D$,
and
$${\Delta(t\theta)\over\Delta_C(Nt)}
=\left|{(\alpha-\sigma\alpha)^2\over{\alpha\sigma\alpha}}\cdot
{\tau(\alpha-\sigma\alpha)^2\over{\tau\alpha\tau\sigma\alpha}}\right|^{1/2}
=\left|{(4b\tau bD)^2\over (a^2-b^2D)(\tau a^2-\tau b^2D)}\right|^{1/2}.$$

\heading Characters\endheading

Denote by $f$ (resp. $f_C$) a complex-valued
compactly-supported smooth (thus locally-constant
since $F$ is nonarchimedean) function on $G$
(resp. $C)$. Fix Haar measures on $G$ and on $C$. 

By a $G$-module $\pi$ (resp. $C$-module $\pi_C$)
we mean an admissible representation ([BZ]) of $G$ (resp.
$C)$ in a complex space. An irreducible $G$-module
$\pi$ is called {\it{$\theta$-invariant}} if it is
equivalent to the $G$-module $^\theta\pi$, defined by
$^\theta\pi(g)=\pi(\theta (g))$. In this case there is
an intertwining operator $A$ on the space of $\pi$ with
$\pi(g)A=A\pi(\theta (g))$ for all $g$. Since $\theta^2=1$
we have $\pi(g)A^2=A^2\pi(g)$ for all $g$, and since
$\pi$ is irreducible $A^2$ is a scalar by Schur's lemma.
We choose $A$ with $A^2=1$. This determines $A$ up to a
sign. When $\pi$ has a Whittaker model, which happens for
all components of cuspidal automorphic representations of 
the adele group $\PGL(4,\A)$, we specify a normalization of 
$A$ which is compatible with a global normalization, as follows,
and then put $\pi(g\times\theta)=\pi(g)\times A$. 

Fix a nontrivial character $\upsi$ of $F$ in $\C^\times$,
and a character $\psi(u)=\upsi(a_{1,2}+a_{2,3}-a_{3,4})$ of
$u=(u_{i,j})$ in the upper triangular subgroup $U$ of $G$.
Note that $\psi(\theta(u))=\psi(u)$. Assume that $\pi$ is a
{\it nondegenerate} $G$-module, namely it embeds in the space
of ``Whittaker'' functions $W$ on $G$, which satisfy -- by
definition -- $W(ugk)=\psi(u)W(g)$ for all $g\in G$, $u\in U$,
$k$ in a compact open subgroup $K_W$ of $K$, as a $G$-module
under right shifts: $(\pi(g)W)(h)=W(hg)$. Then ${}^\theta\pi$ 
is nondegenerate and can be realized in the space of functions
${}^\theta W(g)=W(\theta(g))$, $W$ in the space of $\pi$. We take
$A$ to be the operator on the space of $\pi$ which maps $W$ to
${}^\theta W$.

A $G$-module $\pi$ is called {\it unramified} if the 
space of $\pi$ contains a nonzero $K$-fixed vector. 
The dimension of the space of $K$-fixed vectors is
bounded by one if $\pi$ is irreducible. If $\pi$
is $\theta$-invariant and unramified, and $v_0\neq 0$
is a $K$-fixed vector in the space of $\pi$, then
$Av_0$ is a multiple of $v_0$ (since $\theta K=K$),
namely $Av_0=cv_0$, with $c=\pm 1$. Replace $A$ by
$cA$ to have $Av_0=v_0$, and put $\pi(\theta)=A$.

When $\pi$ is (irreducible) unramified and has a Whittaker 
model, both normalizations of the intertwining operator are
equal. In this case $\upsi$ is unramified (trivial on $R$
but not on $\upi^{-1}R$, where $\upi$ is a generator of the
maximal ideal of $R$), and there exists a unique Whittaker
function $W_0$ in the space of $\pi$ with respect to $\psi$
with $W_0=1$ on $K$. It is mapped by $\pi(\theta)=A$ to ${}^\theta W_0$,
which satisfies ${}^\theta W_0(k)=1$ for all $k$ in $K$ since
$K$ is $\theta$-invariant. Namely $A$ maps the unique normalized
(by $W_0(K)=1$) $K$-fixed vector $W_0$ in the space of $\pi$ to the 
unique normalized $K$-fixed vector ${}^\theta W_0$ in the space of 
${}^\theta\pi$, and we have ${}^\theta W_0=W_0$.

For any $\pi$ and $f$ the convolution operator
$\pi(f)=\int_Gf(g)\pi(g)dg$ has
finite rank. If $\pi$ is $\theta$-invariant put
$\pi(f\times\theta)=\int_Gf(g)\pi(g)\pi(\theta)dg$.
Denote by tr $\pi(f\times\theta)$ the trace of
the operator $\pi(f\times\theta)$. It depends on
the choice of the Haar measure $dg$, but the
({\it{twisted}}) {\it{character}} $\chi_\pi$ of
$\pi$ does not; $\chi_\pi$ is a locally-integrable
(at least in characteristic zero)
complex-valued function on $G\times\theta$ (see [C], [H])
which is $\theta$-conjugacy invariant and locally-constant 
on the $\theta$-regular set, with $\tr\,\pi(f\times\theta)
=\int_Gf(g)\chi_\pi(g\times\theta)dg$ for all $f$.

Local integrability is not used in this work; rather 
it is recovered for our twisted character.

\heading Small representation\endheading

To describe the $G$-module of interest in this paper, note
that a Levi subgroup $M$ of a maximal parabolic subgroup
$P$ of $G$ of type $(3,1)$ is isomorphic to $\GL(3,F)$. 
Hence a $GL(3,F)$-module $\pi_1$ extends to a $P$-module
trivial on the unipotent radical $N (=F^3)$ of $P$. Let
$\delta$ denote (as above) the character $\delta(p)=|\Ad(p)|\Lie N|$ 
of $P$; it is trivial on $N$. Take $P$ to be the upper
triangular parabolic subgroup of type (3,1), and
$M=\{m=\diag(ah,a)^\ast;\,h\in\GL(3,F),\,a\in F^\times\}$.
Here $g^\ast$ denotes the image in $\PGL(4,F)$ of $g$ from 
$\GL(4,F)$. Then the value of $\delta$ at $p=mn$ is 
$\vert\det h\vert$. Denote by $I(\pi_1)$ the $G$-module 
$\pi=\Ind(\delta^{1/2}\pi_1$; $P,G)$ normalizedly induced from 
$\pi_1$ on $P$ to $G$. It is clear from [BZ] that when $\pi_1$ is
self-contragredient and $I(\pi_1)$ is irreducible then it is 
$\theta$-invariant, and it is unramified if and only if $\pi_1$ 
is unramified.

{\it Our aim in this work is to compute the $\theta$-twisted 
character $\chi_\pi$ of the $\PGL(4,F)$-module $\pi=I(1_3)$, 
where $1_3$ is the trivial $P$-module, by purely local means.}

We begin by describing a useful model of our representation,
in analogy with the model of [FK] of an analogous representation 
$I_{(2,1)}(1_2)$ of $\PGL(3,F)$. Indeed we shall express $\pi$
as an integral operator in a convenient model, and integrate
the kernel over the diagonal to compute the character of $\pi$.

Denote by $\mu=\mu_s$ the character
$\mu(x)=\vert x\vert^{(s+1)/2}$ of $F^\times$.
It defines a character $\mu_P=\mu_{s,P}$ of $P$,
trivial on $N$, by $\mu_P(p)=\mu((\det m_3)/m_1{^3})$
if $p=mn$ and $m=\pmatrix m_3&0\\0&m_1\endpmatrix^\ast$
with $m_3$ in $\GL(3,F)$, $m_1$ in $\GL(1,F)$. If $s=0$,
then $\mu_P=\delta^{1/2}$. Let $W_s$ be the space
of complex-valued smooth functions $\psi$ on $G$
with $\psi(pg)=\mu_P(p)\psi(g)$ for all $p$ in $P$
and $g$ in $G$. The group $G$ acts on $W_s$
by right translation: $(\pi_s(g)\psi)(h)=\psi(hg)$.
By definition, $I(1_3)$ is the $G$-module $W_s$
with $s=0$. The parameter $s$ is introduced for
purposes of analytic continuation.

We prefer to work in another model $V_s$ of the
$G$-module $W_s$. Let $V$ denote the space of
column 4-vectors over $F$. Let $V_s$ be the space
of smooth complex-valued functions $\phi$ on
$V-\{0\}$ with $\phi(\lambda\v)=\mu(\lambda)^{-4}\phi(\v)$.
The expression $\mu(\det g)\phi(^tg\v)$, which is
initially defined for $g$ in $\GL(4,F)$, depends only
on the image of $g$ in $G$. The group $G$ acts
on $V_s$ by $(\tau_s(g)\phi)(\v)=\mu(\det g)\phi(^tg\v)$.
Let $\v_0\neq 0$ be a vector of $V$ such that the line
$\{\lambda\v_0;\lambda$ in $F\}$ is fixed under
the action of $^tP$. Explicitly, we take $\v_0=\,{}^t(0,0,0,1)$.
It is clear that the map $V_s\to W_s$, $\phi\mapsto\psi=\psi_\phi$,
where $\psi(g)=(\tau_s(g)\phi)(\v_0)=\mu(\det g)\phi(^tg\v_0)$,
is a $G$-module isomorphism, with inverse $\psi\mapsto\phi=\phi_\psi$,
$\phi(\v)=\mu(\det g)^{-1}\psi(g)$ if $\v=\,{}^tg\v_0$ ($G$
acts transitively on $V-\{0\})$.

For $\v=\,{}^t(x,y,z,t)$ in $V$ put $\Vert\v\Vert=\max(\vert x\vert$,
$\vert y\vert$, $\vert z\vert$, $\vert t\vert)$. Let $V^0$ be the 
quotient of the set of $\v$ in $V$ with $\Vert\v\Vert=1$ by the
equivalence relation $\v\sim\alpha\v$ if $\alpha$ is a
unit in $R$. Denote by $\Bbb P V$ the projective space
of lines in $V-\{0\}$. If $\Phi$ is a function on $V-\{0\}$
with $\Phi(\lambda\v)=\vert\lambda\vert^{-4}\Phi(\v)$ and
$d\v=dx~dy~dz~dt$, then $\Phi(\v)d\v$ is homogeneous of degree
zero. Define
$$\smallint\limits_{\Bbb P V}\Phi(\v)d\v\quad\text{to be}\quad
\smallint\limits_{V^0}\Phi(\v)d\v.$$
Clearly we have
$$\smallint\limits_{\Bbb P V}\Phi(\v)d\v=\smallint\limits_{\Bbb P V}
\Phi(g\v)d(g\v)=\vert\det g\vert\smallint\limits_{\Bbb P V}
\Phi(g\v)d\v.$$
Put $\nu(x)=\vert x\vert$ and $m=2(s-1)$. Note that
$\nu/\mu_s=\mu_{-s}$. Put $\langle \v,\w\rangle =\,{}^t\v J\w$. Then
$\langle g\v,\theta(g)\w\rangle =\langle \v,\w\rangle $.

\proclaim {1. Proposition} The operator $T_s:V_s\to V_{-s}$,
$(T_s\phi)(\v)=\smallint\limits_{\Bbb P V}
\phi(\w)\vert\langle\w,\v\rangle\vert^md\w$,
converges when $\Re s > 1/2$. It satisfies
$T_s\tau_s(g)=\tau_{-s}(\theta(g))T_s$ for all
$g$ in $G$ where it converges. \endproclaim

\demo {Proof} We have
$$\eqalign{(T_s&(\tau_s(g)\phi))(\v)
=\smallint(\tau_s(g)\phi)(\w)\vert^t\w J\v\vert^md\w
=\mu(\det g)\smallint\phi(^tg\w)\vert^t\w J\v\vert^md\w\cr
&=\vert\det g\vert^{-1}\mu(\det g)\smallint\phi(\w)
\vert^t(^tg^{-1}\w)J\v\vert^md\w\cr
&=(\mu/\nu)(\det g)\smallint\phi(\w)\vert^t\w J
\cdot J^{-1}g^{-1}J\v\vert^md\w\cr
&=(\mu/\nu)(\det g)\smallint\phi(\w)\vert\langle\w,
\theta(^tg)\v\rangle\vert^md\w=(\nu/\mu)(\det\theta(g))\cdot
(T_s\phi)({}^t\theta(g)\v)\cr
&=[(\tau_{-s}(\theta (g)))(T_s\phi)](\v)\cr}$$
for the functional equation. 

For the convergence, we may assume that $\phi=1$ and ${}^t\v=(0,0,0,1)$, 
so that the integral is $\int_R\,|x|^mdx$, which converges for $m>-1$. 
Our $m$ is $2s-2$, as required.
\enddemos

The spaces $V_s$ are isomorphic to the space $W$ of
locally-constant complex-valued functions on $V^0$,
and $T_s$ is equivalent to an operator $T^0_s$ on $W$.
The proof of Proposition 1 implies also 

\proclaim {1. Corollary} The operator $T^0_s\circ\tau_s(g^{-1})$
is an integral operator with kernel
$$(\mu/\nu)(\det\theta (g))\vert\langle\w,\theta(^tg^{-1})\v\rangle
\vert^m\qquad (\v,\w~\text{in}~V^0)$$
and trace
$$\tr[T^0_s\circ\tau_s(g^{-1})]=(\nu/\mu)(\det g)
\smallint\limits_{V^0}\vert^t\v gJ\v\vert^md\v.$$
\endproclaim

\remark {Remark} (1)  In the domain where the
integral converges, it is clear that
$\tr[\tau_{-s}({}^tg)\circ T^0_s]$, which is
$\tr[T^0_s\circ\tau_s(g^{-1})]$, depends only on
the $\sigma$-conjugacy class of $g$ if (and only if)
$s=0$.  (2)  To compute the trace of the analytic 
continuation of  $T^0_s\circ\tau_s(g^{-1})$ it 
suffices to compute this trace for $s$ in the domain
of convergence, and then evaluate the resulting
expression at the desired $s$.  Indeed, for each
compact open $\sigma$-invariant subgroup $K$ of
$G$ the space $W_K$ of $K$-biinvariant functions
in $W$ is finite dimensional.  Denote by
$p_K: W\to W_K$ the natural projection.  Then
$p_K\circ T^0_s\circ\tau_s(g^{-1})$ acts on $W_K$,
and the trace of the analytic continuation of
$p_K\circ T^0_s\circ\tau_s(g^{-1})$ is the analytic
continuation of the trace of
$p_K\circ T^0_s\circ\tau_s(g^{-1})$.  Since $K$
can be taken to be arbitrarily small the claim
follows.
\endremark

Next we normalize the operator $T=T_0$ so that it acts trivially on
the one-dimensional space of $K$-fixed vectors in $V_s$. This space
is spanned by the function $\phi_0$ in $V_s$ with $\phi_0(\v)=1$ for
all $\v$ in $V^0$.

Denote again by $\pii$ a generator of the maximal ideal of the ring 
$R$ of integers in our local nonarchimedean field $F$ of odd residual 
characteristic. Denote by $q$ the number of elements of the residue 
field $R/\pii R$ of $R$. Normalize the absolute value by $|\pii|=q^{-1}$, 
and the measures by $\Vol\{|x|\le 1\}=1$. Then $\Vol\{|x|=1\}=1-q^{-1}$,
and the volume of $V^0$ is $(1-q^{-4})/(1-q^{-1})=1+q^{-1}+q^{-2}+q^{-3}$.

\proclaim{2. Proposition} If $\v_0={}^t(0,0,0,1)$, then
$(T\phi_0)(\v_0)=(1-q^{-2(s+1)})/(1-q^{1-2s})\phi_0(\v_0)$. 
When $s=0$, the constant is $-q^{-1}(1+q^{-1})$.
\endproclaim

\demo{Proof} Indeed,
$$(T\phi_0)(\v_0)=\int_{V^0}\phi_0(\v)|{}^t\v J\v_0|^md\v
=\int_{V^0}|x|^m dxdydzdt$$
$$=\left[\int_{||\v||\le 1}-\int_{||\v||<1}\right]
|x|^mdxdydzdt/\int_{|x|=1}dx$$
$$=(1-q^{-m-4})\int_{|x|\le 1}|x|^m dx/\int_{|x|=1}dx
=(1-q^{-2(s+1)})/(1-q^{1-2s}),$$ 
since $m=2(s-1)$ and $\int_{|x|\le 1}|x|^mdx=(1-q^{-m-1})^{-1}\int_{|x|=1}dx$.
\enddemos

\heading Character computation for type I\endheading

For the $\theta$-conjugacy class of type I, represented by 
$g=t\cdot\diag(\r,\s,\s,\r)$, the product 
$${}^t\v gJ\v=(t,z,x,y)
\left(\matrix a_1\r & 0 & 0 & a_2D\r\\ 0 & b_1\s & b_2D\s & 0\\
0 & b_2\s & b_1\s& 0\\ a_2\r & 0 & 0 & a_1\r\endmatrix\right)
\left(\matrix 0 & 0 & 0 & 1\\ 0 & 0 & 1 & 0\\
0 &-1 & 0 & 0\\-1 & 0 & 0 & 0\endmatrix\right)
\left(\matrix t\\ z\\ x\\ y\endmatrix\right)$$
is equal to
$$-t^2a_2D\r-z^2b_2D\s+x^2b_2\s+y^2a_2\r.$$
Note that the trace is a function of $g$ in the projective group, and $\r$ 
and $\s$ range over a set of representatives for $F^\times/N_{E/F}E^{\times}$.

We need to compute
$$({\nu\over{\mu}})(\det g){\Delta(g\theta)\over{\Delta_C(Ng)}}
\int_{V^0}|{}^t\v gJ\v|^md\v$$
$$={|\r\s|^{1-s}|4a_2b_2D|\over{|(a_1^2-a_2^2D)(b_1^2-b_2^2D)|^{s/2}}}
\int_{V^0}|x^2b_2\s+y^2a_2\r-t^2a_2D\r-z^2b_2D\s|^{2(s-1)}dxdydzdt.$$
This is equal to
$$|{\r\over{\s}}|^{-s} |4D\r'| 
|(({a_1\over{b_2}})^2-({a_2\over{b_2}})^2D)(({b_1\over{b_2}})^2-D)|^{-s/2}
\int_{V^0} |x^2-y^2\r'+t^2D\r'-z^2D|^{2(s-1)} dxdydzdt.$$
Here $\r'=-{a_2\over{b_2}}{\r\over\s}$. As $\r$ ranges over 
$F^\times/N_{E/F}E^\times$, we may rename $\r'$ by $\r$.
We get the product of a factor whose value at $s=0$ 
is 1, the factor $|4D\r|$, and the integral 
$$I_s(\r,D)=\int_{V^0}|x^2-\r y^2-D z^2+\r D t^2|^{2(s-1)} dxdydzdt.$$

\remark{Remark} If $F$ is the field $\R$ of real numbers, then its only
algebraic extension is the field $E=\C=F(\sqrt{D})$ of complex numbers,
and only case I occurs. We may take $D=-1$ and $\r$ to range over the
group $\{\pm1\}$. If $\r=-1$ we get 
$$I_s(-1,-1)=\int_{V^0}|x^2+y^2+z^2+t^2|^{2(s-1)}dxdydzdt,$$ 
which is 1 (and has to be multiplied by 4). Is 
$$I_s(1,-1)=\int_{V^0}|x^2-y^2+z^2-t^2|^{2(s-1)}dxdydzdt$$
equal to $-1$ (at least at $s=0$)?
\endremark

\proclaim{I. Theorem} The value of $|4\r D|I_s(\r,D)/(T\phi_0)(\v_0)$ 
at $s=0$ is $2\kappa_E(\r)$, where $\kappa_E$ is the nontrivial
character of $F^\times/N_{E/F}E^\times$, $E=F(\sqrt D)$.
\endproclaim

\demo{Proof} Consider the case when the quadratic form
$x^2-\r y^2-D z^2+\r D t^2$ does not represent zero (is anisotropic). 
Thus $D=\pii$ and $\r\in R^\times-R^{\times 2}$ (hence $|\r|=1$, 
$|D|=1/q$), or $D\in R^\times-R^{\times 2}$ and $\r=\pii$. The second 
case being equivalent to the first, it suffices to deal with the first 
case. The domain 
$\max\{|x|,|y|,|z|,|t|\}=1$ is the disjoint union of
$\{|x|=1\}$, $\{|x|<1,|y|=1\}$, $\{|x|<1,|y|<1,|z|=1\}$ 
and $\{|x|<1,|y|<1,|z|<1,|t|=1\}$.
Thus the integral $I_s(\r,D)$ is the quotient by $\int_{|x|=1}dx$ of
$$\int_{|x|=1} dx+\int\int_{|x|<1,|y|=1} dx dy
+q^{-m}\int\int\int_{|x|<1,|y|<1,|z|=1} dx dy dz$$
$$+q^{-m}\int\int\int\int_{|x|<1,|y|<1,|z|<1,|t|=1}dxdydzdt$$
$$=1+q^{-1}+q^{-m-2}+q^{-m-3}
=1+q^{-1}+q^{-2s}+q^{-2s-1}.$$
The value at $s=0$ is $2(1+q^{-1})$. Since $|\r D|=q^{-1}$, using
Proposition 2 the value of the expression to be evaluated in the
theorem is $-2$. Since $\kappa_E(\r)=-1$, the theorem follows when 
the quadratic form is anisotropic.
\enddemo

We then turn to the case when the quadratic form is isotropic.
Recall that $\r$ ranges over a set of representatives for
$F^\times/N_{E/F}E^\times$, $E=F(\sqrt D)$.
Thus $D\in F-F^2$, and we may assume that $|D|$ and $|\r|$ 
lie in $\{1,q^{-1}\}$. 

\proclaim{I.1. Proposition} When the quadratic form 
$x^2-\r y^2-D z^2+\r D t^2$ is isotropic, $\r$ lies in $N_{E/F}E^\times$,
and we may assume that the quadratic form takes one of three shapes: 
$$x^2-y^2-D z^2+D t^2,\,\,\, D\in R^\times-R^{\times 2};\quad
x^2+\pii y^2-\pii z^2-\pii^2t^2;\quad x^2-y^2-\pii z^2+\pii t^2.$$
\endproclaim

\demo{Proof} 
(1) If $E/F$ is unramified, then $|D|=1$, thus 
$D\in R^\times-R^{\times 2}$. The norm group $N_{E/F}E^\times$ 
is $\pi^{2\Z}R^\times$. If $x^2-\r y^2-D z^2+\r D t^2$ represents 
0 then $\r\in R^\times$, so we may take $\r=1$. 

\n (2) If $E/F$ is ramified then $|D|=q^{-1}$ and 
$N_{E/F}E^\times=(-D)^{\Z}R^{\times 2}$. The form
$x^2-\r y^2-D z^2+\r D t^2$ represents zero when $\r\in 
R^{\times 2}$ or $\r\in -D R^{\times 2}$. Then the
form can be taken to be $x^2+D y^2-D z^2-D^2t^2$ 
with $\r=-D$ and $|D\r|=q^{-2}$, or 
$x^2-y^2-D z^2+D t^2$ with $\r=1$ and
$|D\r|=q^{-1}$. The proposition follows.
\enddemos

The set $V^0=V/\sim$, where 
$V=\{\v=(x,y,z,t)\in R^4;\max\{|x|,|y|,|z|,|t|\}=1\}$
and $\sim$ is the equivalence relation $\v\sim\alpha\v$ for 
$\alpha\in R^\times$, is the disjoint union of the subsets 
$$V_n^0=V_n^0(\r,D)=V_n(\r,D)/\sim,$$ where
$$V_n=V_n(\r,D)=\{\v;\max\{|x|,|y|,|z|,|t|\}=1,
|x^2-\r y^2-D z^2+\r D t^2|=1/q^n\},$$
over $n\ge 0$, and of $\{\v; x^2-\r y^2-D z^2+\r D t^2=0\}/\sim$, 
a set of measure zero.

Thus the integral $I_s(\r,D)$ coincides with the sum
$$\sum_{n=0}^{\infty}q^{-nm}\Vol(V_n^0(\r,D)).$$
When the quadratic form represents zero the problem 
is then to compute the volumes 
$$\Vol(V_n^0(\r,D))=\Vol(V_n(\r,D))/(1-1/q)\qquad (n\ge 0).$$

We need the following Technical Lemma.

\proclaim{I.0. Lemma} When $c\in R^{\times 2}$, $|c|=1$, and $n\ge 1$, we have
$$\int_{|c-x^2|=q^{-n}}dx=\frac{2}{q^n}\left(1-\frac{1}{q}\right),\quad
\text{hence}\quad\int_{|c-x^2|\le q^{-n}}dx=\frac{2}{q^n}.$$
\endproclaim

\demo{Proof} Recall that any $p$-adic number $a$ such that $|a|\le 1$ 
can be written (not uniquely) as a power series in $\pii$:
$$a=\sum_{i=0}^{\infty}a_i\pii^i=a_0+a_1\pii+a_2\pii^2+\dots,\qquad a_i\in R.$$
If $|a|=1/q^n$ we may assume that 
$a_0=a_1=\dots=a_{n-1}=0$ and $a_n\ne 0$. We can write
$$x=\sum_{i=0}^{\infty}x_i\pii^i,\qquad 
c=\sum_{i=0}^{\infty}c_i\pii^i,\qquad 
x^2=\sum_{i=0}^{\infty}a_i\pii^i,\qquad 
a_i=\sum_{j=0}^{i}x_j x_{i-j}\qquad (x_i,c_i,a_i\in R).$$
We have $$c-x^2=\sum_{i=0}^{\infty}(c_i-a_i)\pii^i.$$
Since $|c-x^2|=1/q^n$ we may change the $c_i$ so that
$$c_i-a_i=0\qquad (i=0,...,n-1),\qquad c_n-a_n\ne 0.$$ 

From $c_0=a_0=x_0^2$ it follows that $x_0=\pm c_0'$, where $c_0'$
is a fixed square root of $c_0$ in $R$ (since $c\in R^\times{}^2$,
we have $c_0\in R^\times{}^2$). From
$c_i-a_i=0$ ($i=1,...,n-1$) it follows that (since $x_0\ne 0$)
$$x_i=(c_i-\sum_{j=1}^{i-1}x_j x_{i-j})/(2x_0),\qquad
x_n\ne (c_n-\sum_{j=1}^{n-1}x_j x_{n-j})/(2x_0),$$
where in the case of $i=1$ the sum over $j$ is empty. 
Thus we have
$$\int_{|c-x^2|=q^{-n}}dx=\frac{2}{q}\left(\frac{1}{q}\right)^{n-1}
\left(1-\frac{1}{q}\right)=\frac{2}{q^n}\left(1-\frac{1}{q}\right).$$
The lemma follows.
\enddemos

\proclaim{I.1. Lemma} When $D=\pii$ and $\r=1$, thus $|\r D|=1/q$,
we have
$$\Vol(V_n^0)=\cases
1-1/q,               &\qquad\text{if}\ n=0,\\
q^{-1}(1-1/q)(2+1/q) , &\qquad\text{if}\ n=1,\\
2q^{-n}(1-1/q)(1+1/q), &\qquad\text{if}\ n\ge 2.
\endcases$$
\endproclaim

\demo{Proof} In our case
$$V_0=V_0(1,\pii)=\{(x,y,z,t);\,\max\{|x|,|y|,|z|,|t|\}=1,\, 
|x^2-y^2-\pii z^2+\pii t^2|=1\}.$$
Since $|z|\le 1$ and $|t|\le 1$, we have $|\pii (z^2-t^2)| < 1$, and
$$1=|x^2-y^2-\pii z^2+\pii t^2|=|x^2-y^2|=|x-y||x+y|.$$
Thus $|x-y|=|x+y|=1$, and if $|x|\neq |y|$, $|x\pm y|=\max\{|x|,|y|\}$.
We split $V_0$ into three distinct subsets, corresponding to the cases
$|x|=|y|=1$; $|x|=1$, $|y|<1$; and $|x|<1$, $|y|=1$.
The volume is then
$$\Vol(V_0)=\int_{|t|\le 1}\int_{|z|\le 1}\int_{|x|=1}
\left[\int_{|y|=1,|x-y|=|x+y|=1}\right]dydxdzdt$$
$$+\int_{|t|\le 1}\int_{|z|\le
1}\left[\int_{|x|=1}\int_{|y|<1}+\int_{|x|<1}\int_{|y|=1}
\right]dydxdzdt$$
$$=\int_{|x|=1}\left[\int_{|y|=1,|x-y|=|x+y|=1}\right]dydx
+\frac{2}{q}\left(1-\frac{1}{q}\right)
=\left(1-\frac{1}{q}\right)^2.$$

{\it Let us consider the case of} $V_n$ with $n\ge 2$. If $|x|=1$, then
put $c=c(x,t,z)=x^2+\pii (t^2-z^2)$. Since $|\pii (t^2-z^2)|<1$,
we have $c(x,t,z)\in R^{\times 2}$, and we can apply Lemma I.0. 
Thus we obtain
$$\int_{|t|\le 1}\int_{|z|\le 1}\int_{|x|=1}\int_{|c-y^2|=q^{-n}}
dydxdzdt=\left(1-\frac{1}{q}\right)\frac{2}{q^n}
\left(1-\frac{1}{q}\right)=\frac{2}{q^n}\left(1-\frac{1}{q}\right)^2.$$

If $|x|<1$ it follows that $|y|<1$. Since $\max\{|x|,|y|,|z|,|t|\}=1$
and $n\ge 2$ it follows that $|t|=|z|=1$. Indeed, if, say, $|t|=1$
but $|z|<1$, then $|x^2-y^2-\pii z^2+\pii t^2|=|\pii t^2|=1/q$,
which is a contradiction. Further, dividing by $\pii$ we obtain 
$|z^2+(y^2-x^2)/\pii-t^2|=q^{1-n}$.
Put $c=c(x,y,z)=z^2+(y^2-x^2)/\pii$. Since $|(y^2-x^2)/\pii|<1$, we 
have $c\in R^{\times 2}$, and using Lemma I.0 we obtain
$$\int_{|x|<1}\int_{|y|<1}\int_{|z|=1}\int_{|c-t^2|=q^{1-n}}
dtdzdydx=\frac{1}{q^2}\left(1-\frac{1}{q}\right)\frac{2}{q^{n-1}}
\left(1-\frac{1}{q}\right)=\frac{2}{q^{n+1}}\left(1-\frac{1}{q}\right)^2.$$
Adding the two cases we have ($n\ge 2$)
$$\Vol(V_n)=\frac{2}{q^n}\left(1-\frac{1}{q}\right)^2+
\frac{2}{q^{n+1}}\left(1-\frac{1}{q}\right)^2=
\frac{2}{q^n}\left(1+\frac{1}{q}\right)\left(1-\frac{1}{q}\right)^2.$$

{\it Let us consider the case} $n=1$. The case $|x|=1$, is exactly the
same as for $n\ge 2$. The contribution is $2/q(1-1/q)^2$. Now if
$|x|<1$ then $|y|<1$ and $\max\{|z|,|t|\}=1$. We have
$|x^2-y^2-\pii z^2+\pii t^2|=|\pii(z^2-t^2)|=q^{-1}$. Dividing
by $\pii$ gives $|z^2-t^2|=1$. The volume of this subset is
$$\frac{1}{q^2}\left[\int_{|z|=1}\int_{|z^2-t^2|=1}dtdz+
\int_{|z|<1}\int_{|t|=1}dtdz\right]$$
$$=\frac{1}{q^2}\left[\left(1-\frac{1}{q}\right)\left(1-\frac{2}{q}\right)
+\frac{1}{q}\left(1-\frac{1}{q}\right)\right]=\frac{1}{q^2}
\left(1-\frac{1}{q}\right)^2.$$
Adding the two cases, we have
$$\Vol(V_1)=\frac{2}{q}\left(1-\frac{1}{q}\right)^2+
\frac{1}{q^2}\left(1-\frac{1}{q}\right)^2=
\frac{1}{q}\left(2+\frac{1}{q}\right)\left(1-\frac{1}{q}\right)^2.$$
The lemma follows.
\enddemos

\medskip

\proclaim{I.2. Lemma} When $D=\pii$ and $\r=-\pii$, thus
$|\r D|=1/q^2$, we have
$$\Vol(V_n^0)=\cases
1,                    &\qquad\text{if}\ n=0,\\
q^{-1}(1-1/q),        &\qquad\text{if}\ n=1,\\
q^{-2}(2-1/q-2/q^2),  &\qquad\text{if}\ n=2,\\
2q^{-n}(1-1/q)(1+1/q),&\qquad\text{if}\ n\ge 3.
\endcases $$
\endproclaim

\demo{Proof} {\it To compute} $\Vol(V_0)$, recall that in our case
$$V_0=\{(x,y,z,t);\max\{|x|,|y|,|z|,|t|\}=1, 
|x^2+\pii (y^2-z^2)-\pii^2 t^2|=1\}.$$
Since $|y|\le 1$, $|z|\le 1$, $|t|\le 1$, we have 
$|x^2+\pii (y^2-z^2)-\pii^2 t^2|=|x^2|=1$, and so
$$\Vol(V_0)=\int_{|t|\le 1}\int_{|z|\le 1}\int_{|y|\le 1}\int_{|x|=1}
dxdydzdt=1-\frac{1}{q}.$$

{\it To compute} $\Vol(V_n)$, $n\ge 1$, recall that
$$V_n=\{(x,y,z,t);\max\{|x|,|y|,|z|,|t|\}=1, 
|x^2+\pii (y^2-z^2)-\pii^2 t^2|=1/q^n\}.$$
Assume that $|x|=1$. Then 
$$1=|x^2|=|x^2+\pii (y^2-z^2)-\pii^2 t^2|=1/q^n<1.$$
Thus we have that $|x|<1$ and $\max\{|y|,|z|,|t|\}=1$.

{\it Consider the case} $n=1$. Then $|y^2-z^2|=|y\pm z|=1$. We have
$$\Vol(V_1)=\int_{|t|\le 1}\int_{|x|<1}\left[\int_{|y|=1}
\int_{|y\pm z|=1}+\int_{|y|<1}\int_{|z|=1}\right]dzdydxdt$$
$$=\frac{1}{q}\left[\left(1-\frac{1}{q}\right)\left(1-\frac{2}{q}\right)
+\frac{1}{q}\left(1-\frac{1}{q}\right)\right]
=\frac{1}{q}\left(1-\frac{1}{q}\right)^2.$$

{\it Consider the case} $n\ge 3$. As in the analogous case of
Lemma I.1, we consider the cases of $|y|=1$ and $|y|<1$. 
Adding the two cases we have ($n\ge 3$)
$$\Vol(V_n)=\frac{2}{q^n}\left(1-\frac{1}{q}\right)^2+
\frac{2}{q^{n+1}}\left(1-\frac{1}{q}\right)=
\frac{2}{q^n}\left(1+\frac{1}{q}\right)\left(1-\frac{1}{q}\right)^2.$$

{\it Consider the case} $n=2$. If $|y|=1$ we apply Lemma I.0 (as in the
case $n\ge 3$) and the contribution is $2/q^2(1-1/q)^2$. Now if
$|y|<1$ then $|z|<1$ and $|t|=1$. 
The contribution from this subset is
$$\int_{|y|<1}\int_{|z|<1}\int_{|t|=1}\int_{|t^2-(x/\pi)^2|=1}dxdtdzdy
=\frac{1}{q^2}\left(1-\frac{1}{q}\right)\frac{1}{q}
\left(1-\frac{2}{q}\right).$$
Adding the two cases (see Lemma I.1 for details), we obtain
$$\Vol(V_2)=\frac{2}{q^2}\left(1-\frac{1}{q}\right)^2+\frac{1}{q^3}
\left(1-\frac{1}{q}\right)\left(1-\frac{2}{q}\right)=\frac{1}{q^2}
\left(1-\frac{1}{q}\right)\left(2-\frac{1}{q}-\frac{2}{q^2}\right).$$
The lemma follows.
\enddemos
\medskip

\proclaim{I.3. Lemma} When $E/F$ is unramified, thus $|\r D|=1$, we have
$$\Vol(V_n^0)=\cases
1-1/q^2,                    &\qquad\text{if}\ n=0,\\
q^{-n}(1-1/q)(1+2/q+1/q^2), &\qquad\text{if}\ n\ge 1.
\endcases$$
\endproclaim

\demo{Proof} {\it First we compute} $\Vol(V_0)$. Recall that in our case
$$V_0=\{(x,y,z,t);\max\{|x|,|y|,|z|,|t|\}=1, 
|x^2-y^2-D (z^2-t^2)|=1\}.$$
Since $|x^2-y^2-D (z^2-t^2)|\le\max\{|x|,|y|,|z|,|t|\}$, 
$$V_0=\{(x,y,z,t)\in R^4; |x^2-y^2-D (z^2-t^2)|=1\}.$$

Make the change of variables $x'=x+y$, $y'=x-y$, $z'=z+t$, $t'=z-t$.
Renaming $x'$, $y'$, $z'$, $t'$ as $x$, $y$, $z$, $t$, we obtain
$$V_0=\{(x,y,z,t)\in R^4; |xy-D zt|=1\}.$$

Assume that $|xy|<1$. Since $|xy-D zt|=1$, it follows that 
$|zt|=|z|=|t|=1$. The contribution from the set $|xy|<1$ is
$$\int_{|t|=1}\int_{|z|=1}\left[\int_{|x|<1}\int_{|y|\le1}+\int_{|x|=1}
\int_{|y|<1}\right]dydxdzdt$$
$$=\left(1-\frac{1}{q}\right)^2\left(\frac{1}{q}+
\left(1-\frac{1}{q}\right)\frac{1}{q}\right)
=\frac{1}{q}\left(1-\frac{1}{q}\right)^2\left(2-\frac{1}{q}\right).$$

Note that the contribution from $|xy|=1$, $|zt|<1$, is the same and
equals
$$\frac{1}{q}\left(1-\frac{1}{q}\right)^2\left(2-\frac{1}{q}\right).$$

We are left with the case $|xy|=|zt|=1$, i.e. $|x|=|y|=|z|=|t|=1$.
If $|x|=|y|=|z|=1$ we introduce $U(x,y,z)=\{t;\ |t|=1,\ |xy-D zt|=1\}$,
a set of volume $1-2/q$. The contribution from this case is
$$\int_{|x|=1}\int_{|y|=1}\int_{|z|=1}\int_{U(x,y,z)}dtdzdydx
=\left(1-\frac{1}{q}\right)^3\left(1-\frac{2}{q}\right).$$
Thus we obtain
$$\Vol(V_0)=\frac{2}{q}\left(1-\frac{1}{q}\right)^2\left(2-\frac{1}{q}\right)
+\left(1-\frac{1}{q}\right)^3\left(1-\frac{2}{q}\right)=
\left(1-\frac{1}{q}\right)^2\left(1+\frac{1}{q}\right).$$

{\it Next we compute} $\Vol(V_n)$, $n\ge 1$. Recall that in our case
$$V_n=\{(x,y,z,t);\max\{|x|,|y|,|z|,|t|\}=1, 
|x^2-y^2-D (z^2-t^2)|=1/q^n\}.$$

Making the change of variables $u=x+y$, $v=x-y$, we obtain
$$V_n=\{(u,v,z,t);\max\{|u+v|,|u-v|,|z|,|t|\}=1,
|uv-D (z^2-t^2)|=1/q^n\}.$$
Since the set $\{v=0\}$ is of measure zero, we assume that
$v\ne 0$. Then $|uv-D (z^2-t^2)|=1/q^n$ implies that 
$u=D (z^2-t^2) v^{-1}+wv^{-1}\pii^n$, where $|w|=1$. 
There are two cases.

Assume that $|v|=1$. Note that if
$|z^2-t^2|=1$, then $\max\{|z|,|t|\}=1$, and if
$|z^2-t^2|<1$, then (since $n\ge 1$)
$$|u|=|D (z^2-t^2) v^{-1}+wv^{-1}\pii^n|\le\max\{|z^2-t^2|,q^{-n}\}<1,$$
and consequently $|u+v|=|v|=1$. So $|v|=1$ implies that
$\max\{|u+v|,|u-v|,|z|,|t|\}=1$. Further, since $|v|=1$, we have 
$du=q^{-n}dw$. Thus the contribution from the set with $|v|=1$ is
$$\int_{|t|\le 1}\int_{|z|\le 1}\int_{|v|=1}
\int_{|uv-D (z^2-t^2)|=1/q^n}dudvdzdt
=\int_{|v|=1}\int_{|w|=1}\frac{dw}{q^n}dv
=\frac{1}{q^n}\left(1-\frac{1}{q}\right)^2.$$

Assume that $|v|<1$ and $|u|=1$. Thus $\max\{|u+v|,|u-v|,|z|,|t|\}=1$.
We write $v=D (z^2-t^2) u^{-1}+wu^{-1}\pii^n$ where $|w|=1$
and $dv=q^{-n}dw$. Since 
$$|z^2-t^2|\le\max\{|v|,q^{-n}\}< 1,$$
it follows that $|z^2-t^2|<1$. Note that
$$\int\int_{|z^2-t^2|<1}dzdt=\int_{|z|=1}\int_{|z^2-t^2|<1}dtdz
+\int_{|z|<1}\int_{|t|<1}dt=\left(1-\frac{1}{q}\right)\frac{2}{q}
+\frac{1}{q^2}=\frac{1}{q}\left(2-\frac{1}{q}\right).$$
The volume of this subset equals
$$\int\int_{|z^2-t^2|<1}\int_{|u|=1}\int_{|uv-D (z^2-t^2)|=1/q^n}
dvdudzdt=\frac{1}{q}\left(2-\frac{1}{q}\right)\frac{1}{q^n}
\int_{|u|=1}\int_{|w|=1}dwdu$$
$$=\frac{1}{q}\frac{1}{q^n}\left(1-\frac{1}{q}\right)^2
\left(2-\frac{1}{q}\right).$$

Assume that $|v|<1$ and $|u|<1$. Then we have $|u\pm v|<1$ and
thus $\max\{|z|,|t|\}=1$. Since $|uv-D (z^2-t^2)|<1$ it 
follows that $|z^2-t^2|<1$. So we have $|z|=|t|=1$. Put 
$c=c(z,u,v)=z^2-uvD^{-1}$. Then $c\in R^{\times 2}$ 
(since $|uvD^{-1}|<1$). Dividing by $D$, we have
$$\frac{1}{q^n}=|uv-D (z^2-t^2)|=|c-t^2|.$$
Applying Lemma I.0, the contribution from this subset is equal to
$$\int_{|u|<1}\int_{|v|<1}\int_{|z|=1}\int_{|c-t^2|=q^{-n}}dtdzdvdu
=\frac{1}{q^2}\frac{2}{q^n}\left(1-\frac{1}{q}\right)^2.$$
 
Adding the contributions from $|v|=1$ and $|v|<1$ we obtain
$$\Vol(V_n)=\frac{1}{q^n}\left(1-\frac{1}{q}\right)^2
+\frac{1}{q}\frac{1}{q^n}\left(1-\frac{1}{q}\right)^2
\left(2-\frac{1}{q}\right)
+\frac{1}{q^2}\frac{2}{q^n}\left(1-\frac{1}{q}\right)^2$$
$$=\frac{1}{q^n}\left(1-\frac{1}{q}\right)^2
\left(1+\frac{2}{q}+\frac{1}{q^2}\right).$$
The lemma follows.
\enddemos
\bigskip

\demo{Proof of Theorem I} We are now ready to complete the
proof of Theorem I in the isotropic case. Recall that
we need to compute the value at $s=0$ ($m=-2$) of the product 
$|\r D|I_s(\r,D)$. Here $I_s(\r,D)$ coincides with the sum
$$\sum_{n=0}^{\infty}q^{-nm}\Vol(V_n^0(\r,D))$$
which converges for $m>-1$ by Proposition 1 or alternatively
by Lemmas I.1-I.3. The value at $m=-2$ is obtained then by 
analytic continuation of this sum.

\n {\it Case of Lemma} I.1. We have $|\r D|=1/q$, and 
$I_s(\r,D)$ is equal to
$$\Vol(V_0^0)+q^{-m}\Vol(V_1^0)+\sum_{n=2}^{\infty}q^{-nm}\Vol(V_n^0)$$
$$=1-\frac{1}{q}+\frac{1}{q}\left(1-\frac{1}{q}\right)
\left(2+\frac{1}{q}\right) 
\frac{1}{q^m}+2\left(1-\frac{1}{q}\right)\left(1+\frac{1}{q}\right)
q^{-2(m+1)}\left(1-\frac{1}{q^{m+1}}\right)^{-1}.$$
When $m=-2$, this is 
$$1-\frac{1}{q}+q\left(2-\frac{1}{q}-\frac{1}{q^2}\right) 
+2\left(1-\frac{1}{q}\right)\left(1+\frac{1}{q}\right)\frac{q^2}{1-q}
=-2(1+q^{-1}).$$
Multiplied by $|\r D|=1/q$, we obtain $-2q^{-1}(1+q^{-1})$.

\n {\it Case of Lemma} I.2. We have $|\r D|=1/q^2$, 
and $I_s(\r,D)$ is equal to
$$\Vol(V_0^0)+q^{-m}\Vol(V_1^0)+q^{-2m}\Vol(V_2^0)
+\sum_{n=3}^{\infty}q^{-nm}\Vol(V_n^0)$$
$$=1+\frac{1}{q}\left(1-\frac{1}{q}\right)q^{-m} 
+\frac{1}{q^2}\left(2-\frac{1}{q}-\frac{2}{q^2}\right)q^{-2m}$$
$$+2\left(1-\frac{1}{q}\right)\left(1+\frac{1}{q}\right)
q^{-3(m+1)}\left(1-\frac{1}{q^{m+1}}\right)^{-1}.$$
When $m=-2$, this is 
$$1+\frac{1}{q}\left(1-\frac{1}{q}\right)q^2 
+\frac{1}{q^2}\left(2-\frac{1}{q}-\frac{2}{q^2}\right)q^4
+ 2\left(1-\frac{1}{q}\right)\left(1+\frac{1}{q}\right)
\frac{q^3}{1-q}.$$
Once simplified and multiplied by $|\r D|=1/q^2$, we obtain 
$-2q^{-1}(1+q^{-1})$.

\n {\it Case of Lemma} I.3. We have $|\r D|=1$, and $I_s(\r,D)$ is
equal to
$$\Vol(V_0^0)+\sum_{n=1}^{\infty}q^{-nm}\Vol(V_n^0)$$
$$=1-\frac{1}{q^2}+\left(1-\frac{1}{q}\right)
\left(1+\frac{2}{q}+\frac{1}{q^2}\right)
q^{-(m+1)}\left(1-\frac{1}{q^{m+1}}\right)^{-1}.$$
When $m=-2$, this is 
$$=1-\frac{1}{q^2}+\left(1-\frac{1}{q}\right)
\left(1+\frac{2}{q}+\frac{1}{q^2}\right)\frac{q}{1-q}
=-2q^{-1}(1+q^{-1}).$$

In all 3 cases the value of the expression of the theorem is 2
by virtue of Proposition 2, and indeed $\kappa_E(\r)=1$ in these 
cases as $\r$ lies in $N_{E/F}E^\times$ by Proposition I.1. 
The theorem follows.
\enddemos

\heading Character computation for type II\endheading

For the $\theta$-conjugacy class of type II, represented by 
$g=t\cdot\diag(\r,\s,\s,\r)$, the product 
$${}^t\v gJ\v=(t,z,x,y)
\left(\matrix a_1\r & 0 & 0 & a_2D\r\\ 0 & b_1\s & b_2AD\s & 0\\
0 & b_2\s & b_1\s& 0\\ a_2\r & 0 & 0 & a_1\r\endmatrix\right)
\left(\matrix 0 & 0 & 0 & 1\\ 0 & 0 & 1 & 0\\
0 &-1 & 0 & 0\\-1 & 0 & 0 & 0\endmatrix\right)
\left(\matrix t\\ z\\ x\\ y\endmatrix\right)$$
is equal to
$$-t^2a_2D\r-z^2b_2AD\s+x^2b_2\s+y^2a_2\r,$$
where $a_1+a_2\sqrt D\in E_1^\times$ ($E_1=F(\sqrt{D})$) and
$b_1+b_2\sqrt{AD}\in E_2^\times$ ($E_2=F(\sqrt{AD})$). The trace 
being a function of $g$ in the projective group, and $\r$ ranging 
over a set of representatives for $F^\times/N_{E_1/F}E_1^{\times}$
(and $\s$ for $F^\times/N_{E_2/F}E_2^{\times}$),
we may divide the quadratic form by $b_2\s$ and rename $-a_2\r/b_2\s$
by $\r$. The quadratic form becomes $x^2-y^2\r-z^2AD+t^2D\r$.

Thus we need to compute the integral
$$I_s(\r,A,D)=\int_{V^0}|x^2-\r y^2-AD z^2+\r D t^2|^{2(s-1)} dxdydzdt.$$

The property of the numbers $A$, $D$ and $AD$ that we need is that
their square roots generate the three distinct quadratic extensions
of $F$. Thus we may assume that $\{A,D,AD\}=\{u,\pii,u\pii\}$, where 
$u\in R^\times-R^{\times 2}$. Of course with this normalization $AD$
is no longer the product of $A$ and $D$, but its representative in the
set $\{1,u,\pii,u\pii\}$ mod $F^{\times 2}$. Since $\r$ ranges over a 
set of representatives for $F^\times/N_{E_1/F}E_1^\times$, it can be 
assumed to range over $\{1,\pii\}$ if $D=u$, and over $\{1,u\}$ 
if $|D|=|\pii|$. 

In this section we prove

\proclaim{II. Theorem} The value of $I_s(\r,A,D)$ at $s=0$ is $0$.
\endproclaim

To prove this theorem we need some lemmas.

\proclaim{II.1. Proposition} The quadratic form 
$x^2-\r y^2-AD z^2+\r D t^2$ takes one of six forms: 
$x^2-y^2+\pii(t^2-uz^2)$, 
$x^2-u y^2+u\pii(t^2-z^2)$, 
$x^2-y^2+ut^2-u\pii z^2$, 
$x^2-y^2-u z^2+\pii t^2$,
$x^2-\pii y^2+u\pii(t^2-z^2)$,
$x^2-u y^2-u z^2+u\pii t^2$, 
where $u\in R^\times-R^{\times 2}$.
It is always isotropic.
\endproclaim

\demo{Proof} 
(1) If $E_1/F$ is unramified then $D=u$ where 
$u\in R^\times-R^{\times 2}$. The norm group $N_{E_1/F}E_1^\times$ 
is $\pii^{2\Z}R^\times$. So $\r=1$ or $\pii$ and $A=\pii$.
We obtain two quadratic forms:
$x^2-y^2-u\pii z^2+ut^2$, and $x^2-\pii y^2-u\pii(z^2-t^2)$.

\n (2) If $E_1/F$ is ramified then $D=\pii$ and 
$N_{E_1/F}E_1^\times=(-D)^{\Z}R^{\times 2}$. Then 
$\r=1$ or $u$, and $A=u$ or $u\pii$. Note that if $A=u\pii$
we take $AD=u$. We obtain the following quadratic forms:
if $\r=1$, $A=u$ we have $x^2-y^2-\pii(u z^2-t^2)$;
if $\r=1$, $A=u\pii$ we have $x^2-y^2-u z^2+\pii t^2$;
if $\r=u$, $A=u$ we have $x^2-u y^2-u\pii (z^2-t^2)$;
if $\r=u$, $A=u\pii$ we have $x^2-u y^2-u z^2+u\pii t^2$.
The proposition follows.
\enddemos

The set $V^0=V/\sim$, where 
$V=\{\v=(x,y,z,t)\in R^4;\max\{|x|,|y|,|z|,|t|\}=1\}$
and $\sim$ is the equivalence relation $\v\sim\alpha\v$ for 
$\alpha\in R^\times$, is the disjoint union of the subsets 
$$V_n^0=V_n^0(\r,A,D)=V_n(\r,A,D)/\sim,$$ where
$$V_n=V_n(\r,A,D)=\{\v;\max\{|x|,|y|,|z|,|t|\}=1,
|x^2-\r y^2-AD z^2+\r D t^2|=1/q^n\},$$
over $n\ge 0$, and of $\{\v; x^2-\r y^2-AD z^2+\r D t^2=0\}/\sim$, 
a set of measure zero.

Thus the integral $I_s(\r,A,D)$ coincides with the sum
$$\sum_{n=0}^{\infty}q^{-nm}\Vol(V_n^0(\r,A,D)).$$
The problem is then to compute the volumes 
$$\Vol(V_n^0(\r,A,D))=\Vol(V_n(\r,A,D))/(1-1/q)\qquad (n\ge 0).$$

In the following lemmas, $u$ is a nonsquare unit.  

\proclaim{II.1. Lemma} When the quadratic form is 
$x^2-y^2+\pii (t^2-uz^2)$, we have
$$\Vol(V_n^0)=\cases
1-1/q,           &\qquad\text{if}\ n=0,\\
2/q-1/q^2+1/q^3, &\qquad\text{if}\ n=1,\\
2q^{-n}(1-1/q),  &\qquad\text{if}\ n\ge 2.
\endcases$$
\endproclaim

\demo{Proof} In our case
$$V_0=\{(x,y,z,t);\,\max\{|x|,|y|,|z|,|t|\}=1,\, 
|x^2-y^2+\pii(t^2-uz^2)|=1\}.$$
This is the same case as that of Lemma I.1.
The volume is then
$$\Vol(V_0)=\left(1-\frac{1}{q}\right)^2.$$

{\it Let us consider the case of} $V_n$ with $n\ge 2$. If $|x|=1$, then
put $c=c(x,t,z)=x^2+\pii (t^2-uz^2)$. Since $|\pii (t^2-uz^2)|<1$,
we have $c(x,t,z)\in R^{\times 2}$, and we can apply Lemma I.0. 
Thus we obtain
$$\int_{|t|\le 1}\int_{|z|\le 1}\int_{|x|=1}\int_{|c-y^2|=q^{-n}}
dydxdzdt=\left(1-\frac{1}{q}\right)\frac{2}{q^n}
\left(1-\frac{1}{q}\right)=\frac{2}{q^n}\left(1-\frac{1}{q}\right)^2.$$

If $|x|<1$ it follows that $|y|<1$. Thus $\max\{|z|,|t|\}=1$.
Since $t^2-uz^2$ does not represent zero non trivially, we
have $|t^2-uz^2|=1$, which is a contradiction. Thus
$$\Vol(V_n)=\frac{2}{q^n}\left(1-\frac{1}{q}\right)^2.$$

{\it Let us consider the case} $n=1$. Note that in this case $|x^2-y^2|<1$.

\n (A) Case of $|t^2-uz^2|<1$. Since $t^2-uz^2$ does not represent
zero non trivially, it follows that $|z|<1$, $|t|<1$. Thus
$|\pii(t^2-uz^2)|<1/q^2$, and 
$$V_1=\{(x,y,z,t);\,\max\{|x|,|y|\}=1,\,|x^2-y^2|=1/q\}.$$
Applying Lemma I.0, the contribution of this case is
$$\int_{|z|<1}\int_{|t|<1}\int_{|x|=1}\int_{|x^2-y^2|=q^{-1}}dydxdtdz
=\frac{1}{q^2}\left(1-\frac{1}{q}\right)
\frac{2}{q}\left(1-\frac{1}{q}\right).$$

\n (B1) Case of $|t^2-uz^2|=1$ and $|x^2-y^2|< 1/q$. Since
$t^2-uz^2$ does not represent zero non trivially, we have
$$\int\int_{|t^2-uz^2|=1}dzdt=\int_{|t|=1}\int_{|z|\le 1}dzdt
+\int_{|t|<1}\int_{|z|=1}dzdt=\left(1-\frac{1}{q}\right)
\left(1+\frac{1}{q}\right).$$ 
Furthermore,
$$\int\int_{|x^2-y^2|\le 1/q^2}dxdy=\int_{|x|=1}\int_{|x^2-y^2|
\le1/q^2}dydx+\int_{|x|<1}\int_{|y|<1}dydx=
\left(1-\frac{1}{q}\right)\frac{2}{q^2}+\frac{1}{q^2}.$$
Thus the contribution of this case is equal to
$$\left(1-\frac{1}{q}\right)\left(1+\frac{1}{q}\right) 
\left[\left(1-\frac{1}{q}\right)\frac{2}{q^2}+\frac{1}{q^2}\right].$$

\n (B2) Case of $|t^2-uz^2|=1$ and $|x^2-y^2|=1/q$. Set $w=x-y$,
$v=x+y$ ($dwdv=dxdy$). Then $|wv|=1/q$ and $|wv+\pii(t^2-uz^2)|=1/q$,
and the contribution of this case is (there are two integrals
that correspond to $|w|=1$, $|v|=1/q$ and $|w|=1/q$, $|v|=1$):
$$2\int\int_{|t^2-uz^2|=1}\int_{|v|=1/q}\int_{|w|=1,
|wv+\pi (t^2-uz^2)|=1/q}dwdvdtdz.$$
Since $|(t^2-uz^2)(v/\pii)^{-1}|=1$, we have that 
$w\neq 0,-(t^2-uz^2)(v/\pii)^{-1}$ (mod $\pii$). So, the above
integral is equal to
$$2\left(1-\frac{1}{q}\right)\left(1+\frac{1}{q}\right)\frac{1}{q}
\left(1-\frac{1}{q}\right)\left(1-\frac{2}{q}\right).$$

Adding the contributions from Cases (A), (B1), and (B2) (divided
by $(1-1/q)$), we obtain
$$\Vol(V_1^0)=\frac{2}{q^3}\left(1-\frac{1}{q}\right)+\frac{1}{q^2}
\left(1+\frac{1}{q}\right)\left(3-\frac{2}{q}\right)+\frac{2}{q}
\left(1-\frac{1}{q^2}\right)\left(1-\frac{2}{q}\right).$$
Once simplified this is equal to $2/q-1/q^2+1/q^3$.
The lemma follows.
\enddemos
\medskip

\proclaim{II.2. Lemma} When the quadratic form is $x^2-uy^2+u\pii (t^2-z^2)$, 
we have
$$\Vol(V_n^0)=\cases
1+1/q,             &\qquad\text{if}\ n=0,\\
q^{-2}(1-1/q),     &\qquad\text{if}\ n=1,\\
2q^{-(n+1)}(1-1/q),&\qquad\text{if}\ n\ge 2.
\endcases$$
\endproclaim

\demo{Proof} {\it Consider the case of} $n=0$. Then
$$V_0=\{(x,y,z,t);\,\max\{|x|,|y|,|z|,|t|\}=1,\, 
|x^2-uy^2+u\pii (t^2-z^2)|=1\}.$$
Since $x^2-uy^2$ does not represent zero non trivially, we have
$$1=|x^2-uy^2+u\pii (t^2-z^2)|=|x^2-uy^2|=\max\{|x|,|y|\}.$$
We obtain 
$$\Vol(V_0)=\int_{|x|=1}\int_{|y|\le 1}dydx+\int_{|x|<1}\int_{|y|=1}dydx
=1-\frac{1}{q}+\frac{1}{q}\left(1-\frac{1}{q}\right)
=\left(1-\frac{1}{q}\right)\left(1+\frac{1}{q}\right).$$

{\it Let us consider the case of} $n=1$. Then
$$V_1=\{(x,y,z,t);\,\max\{|x|,|y|,|z|,|t|\}=1,\, 
|x^2-uy^2+u\pii (t^2-z^2)|=1/q\}.$$
If $\max\{|x|,|y|\}=1$, then $|x^2-uy^2|=1$. It
implies that $|x|<1$, $|y|<1$, and $|x^2-uy^2|\le 1/q^2$.
The contribution from this subset is equal to
$$\int_{|x|<1}\int_{|y|<1}\int\int_{|t^2-z^2|=1}dtdz dydx
=\frac{1}{q^2}\left[\int_{|t|=1}\int_{|t^2-z^2|=1} dzdt+
\int_{|t|<1}\int_{|z|=1}dzdt\right]$$
$$=\frac{1}{q^2}\left[\left(1-\frac{1}{q}\right)\left(1-\frac{2}{q}\right)
+\frac{1}{q}\left(1-\frac{1}{q}\right)\right]=\frac{1}{q^2}
\left(1-\frac{1}{q}\right)^2.$$

{\it Let us consider the case of} $n\ge 2$. Then
$$V_n=\{(x,y,z,t);\,\max\{|x|,|y|,|z|,|t|\}=1,\, 
|x^2-uy^2+u\pii (t^2-z^2)|=q^{-n}\}.$$
If $\max\{|x|,|y|\}=1$ then $|x^2-uy^2|=1$, which is
a contradiction. Hence $|x|<1$, $|y|<1$, and 
$\max\{|z|,|t|\}=1$. The latter implies that $|z|=|t|=1$.
Dividing by $u\pii$, we have
$$V_n=\{(x,y,z,t);\,|x|<1,\,|y|<1,\,|z|=1,\,|t|=1,\, 
|z^2-t^2+\pii ((y/\pii)^2-u^{-1}(x/\pii)^2)|=q^{1-n}\}.$$
Applying Lemma I.0 (with $c=z^2+\pii((y/\pii)^2-u^{-1}(x/\pii)^2)$), 
its volume is equal to
$$\int_{|x|<1}\int_{|y|<1}\int_{|z|=1}\int_{|c-t^2|=q^{1-n}}dt dzdydx
=\frac{1}{q^2}\frac{2}{q^{n-1}}\left(1-\frac{1}{q}\right)^2.$$
Dividing by $(1-1/q)$ we obtain the $\Vol(V_n^0)$.
The lemma follows.
\enddemos
\medskip

In all other cases we get the same result for the volumes.
Since the proofs are different, we state the remaining cases
separately as Lemmas II.3, II.4, II.5.

\proclaim{II.3. Lemma} When the quadratic form is $x^2-y^2+ut^2-u\pii z^2$
or $x^2-y^2-uz^2+\pii t^2$, we have
$$\Vol(V_n^0)=\cases
1,              &\qquad\text{if}\ n=0,\\
1/q,            &\qquad\text{if}\ n=1,\\
q^{-n}(1-1/q^2),&\qquad\text{if}\ n\ge 2.
\endcases$$
\endproclaim

\demo{Proof} Since the quadratic form $x^2-y^2-uz^2+\pii t^2$
is equal to $-(y^2-x^2+uz^2-\pii t^2)$, the computations
for this form are identical to those of $x^2-y^2+ut^2-u\pii z^2$.  

{\it Consider the case of} $n=0$. Then
$$V_0=\{(x,y,z,t);\,\max\{|x|,|y|,|z|,|t|\}=1,\, 
|x^2-y^2+ut^2-u\pii z^2|=1\}.$$
We have the following three cases.
\smallskip
\n (A) Case of $|t|<1$. It follows that $|x^2-y^2|=1$. The
contribution of this case is
$$\int_{|z|\le 1}\int_{|t|<1}\int\int_{|x^2-y^2|=1}dydx dtdz
=\frac{1}{q}\left[\int_{|x|=1}\int_{|x^2-y^2|=1}dydx+
\int_{|x|<1}\int_{|y|=1}dydx\right]$$
$$=\frac{1}{q}\left[\left(1-\frac{1}{q}\right)\left(1-\frac{2}{q}\right)
+\frac{1}{q}\left(1-\frac{1}{q}\right)\right]=\frac{1}{q}
\left(1-\frac{1}{q}\right)^2.$$

\n (B1) Case of $|t|=1$, $|x^2-y^2|<1$. The contribution from
this case is equal to (we apply Lemma I.0):
$$\int_{|t|=1}\int_{|x|=1}\int_{|x^2-y^2|\le q^{-1}}dydxdt
+\int_{|t|=1}\int_{|x|<1}\int_{|y|<1}dydxdt$$
$$=\frac{2}{q}\left(1-\frac{1}{q}\right)^2+\frac{1}{q^2}
\left(1-\frac{1}{q}\right)=\left(1-\frac{1}{q}\right)
\left(\frac{2}{q}-\frac{1}{q^2}\right).$$

\n (B2) Case of $|t|=1$, $|x^2-y^2|=1$. Set $w=x-y$, $v=x+y$.
Then $|w|=|v|=1$ and also $|ut^2w^{-1}|=1$. Thus the contribution 
from this case is given by the integral
$$\int_{|t|=1}\int_{|w|=1}\int_{|v|=1,|wv+ut^2|=1} dvdwdt=
\left(1-\frac{1}{q}\right)^2\left(1-\frac{2}{q}\right).$$

Adding the contributions from Cases (A), (B1), and (B2) (divided
by $(1-1/q)$), we obtain
$$\Vol(V_0^0)=\frac{1}{q}\left(1-\frac{1}{q}\right)+\frac{2}{q}
-\frac{1}{q^2}+\left(1-\frac{1}{q}\right)\left(1-\frac{2}{q}\right)=1.$$

{\it Let us consider the case of} $n\ge 2$. We have the following 
three cases.
\smallskip
\n (A) Case of $|x|<1$. Since the quadratic form $y^2-ut^2$
does not represent zero non trivially, we have that $|y^2-ut^2|=1$
if and only if $\max\{|y|,|t|\}=1$. It implies that $|y|<1$, $|t|<1$,
and thus $|x^2-y^2+ut^2-u\pii z^2|=|\pii z^2|=1/q$, which is a
contradiction.
\smallskip
\n (B1) Case of $|x|=1$, $|t|<1$. The contribution from this
case is given by the following integral (we apply Lemma I.0):
$$\int_{|t|<1}\int_{|z|\le 1}\int_{|x|=1}\left(
\int_{|(x^2+ut^2-u\upi z^2)-y^2|=q^{-n}}dy\right)
dxdzdt=\frac{1}{q}\left(1-\frac{1}{q}\right)
\frac{2}{q^n}\left(1-\frac{1}{q}\right).$$

\n (B2) Case of $|x|=1$, $|t|=1$. Set $w=x-y$, $v=x+y$.
Then we have that $|w|=|v|=1$, and from 
$|wv+ut^2-u\pii z^2|=q^{-n}$, we have
$$w=u(\pii z^2-t^2)v^{-1}+\varepsilon v^{-1}\pii^n, 
\qquad dw=\frac{1}{q^n}d\varepsilon,\qquad |\varepsilon|=1.$$
The volume of this subset is given by
$$\int_{|z|\le 1}\int_{|t|=1}\int_{|v|=1}
\int_{|w|=1,|vw+ut^2-u\pi z^2|=q^{-n}}dw dv dtdz$$
$$=\left(1-\frac{1}{q}\right)^2\int_{|\varepsilon|=1}\frac{d\varepsilon}{q^n}
=\frac{1}{q^n}\left(1-\frac{1}{q}\right)^3.$$

Adding the contributions from cases (A), (B1), and (B2) (divided
by $(1-1/q)$), we obtain
$$\Vol(V_n^0)=\frac{2}{q}\left(1-\frac{1}{q}\right)\frac{1}{q^n}
+\left(1-\frac{1}{q}\right)^2\frac{1}{q^n}=\left(1-\frac{1}{q^2}\right)
\frac{1}{q^n}.$$

{\it Let us consider the case of} $n=1$. We have the following
three cases.
\smallskip
\n (A) Case of $|x|<1$. Since the quadratic form $y^2-ut^2$
does not represent zero non trivially, we have that $|y^2-ut^2|=1$
if and only if $\max\{|y|,|t|\}=1$. Hence $|y|<1$, $|t|<1$,
and thus $|z|=1$. The volume of this subset is equal to
$$\int_{|x|<1}\int_{|y|<1}\int_{|t|<1}\int_{|z|=1}dzdtdydx
=\frac{1}{q^3}\left(1-\frac{1}{q}\right).$$

\n (B1) Case of $|x|=1$, $|t|<1$. Applying Lemma I.0, we have
$$\int_{|t|<1}\int_{|z|\le 1}\int_{|x|=1}\left(
\int_{|(x^2+ut^2-u\upi z^2)-y^2|=q^{-1}}dy\right)
dxdzdt=\frac{1}{q}\left(1-\frac{1}{q}\right)\frac{2}{q}
\left(1-\frac{1}{q}\right).$$

\n (B2) Case of $|x|=1$, $|t|=1$. Set $w=x-y$, $v=x+y$. We
have that $|w|=|v|=1$, and we arrive to the same case as that
of $n\ge 2$ (with $n=1$). The contribution is
$$\frac{1}{q}\left(1-\frac{1}{q}\right)^3.$$

Adding the contributions from cases (A), (B1), and (B2) (divided
by $(1-1/q)$), we obtain
$$\Vol(V_1^0)=\frac{1}{q^3}+\frac{2}{q^2}\left(1-\frac{1}{q}\right)
+\frac{1}{q}\left(1-\frac{1}{q}\right)^2=\frac{1}{q}.$$
The lemma follows.
\enddemos
\medskip

\proclaim{II.4. Lemma} When the quadratic form is 
$x^2-\pii y^2+u\pii (t^2-z^2)$, we have
$$\Vol(V_n^0)=\cases
1,              &\qquad\text{if}\ n=0,\\
1/q,            &\qquad\text{if}\ n=1,\\
q^{-n}(1-1/q^2),&\qquad\text{if}\ n\ge 2.
\endcases$$
\endproclaim

\demo{Proof} {\it Consider the case of} $n=0$. Then
$$V_0=\{(x,y,z,t);\,\max\{|x|,|y|,|z|,|t|\}=1,\, 
|x^2-\pii y^2+u\pii (z^2-t^2)|=1\}.$$
Obviously we have 
$$\Vol(V_0)=\int_{|x|=1}dx=1-\frac{1}{q}.$$

{\it Let us consider the case of} $n\ge 2$. It follows that $|x|<1$,
and dividing by $\pii$, we have
$$V_n=\{(x,y,z,t);\,\max\{|y|,|z|,|t|\}=1,\,|x|<1,\, 
|z^2-t^2+uy^2-u\pii (x/\pii)^2|=q^{1-n}\}.$$
This case is the same as that of Lemma II.3. 
We have that $\Vol(V_n^0)$ is the product of $1/q$
and the $\Vol(V_{n-1}^0)$ of Lemma II.3, which is
equal to $q^{-1}(1-1/q^2)q^{-(n-1)}=(1-1/q^2)q^{-n}$.

{\it Let us consider the case of} $n=1$. Then
$$V_1=\{(x,y,z,t);\,\max\{|x|,|y|,|z|,|t|\}=1,\, 
|x^2-\pii y^2+u\pii (z^2-t^2)|=1/q\}.$$
It follows that $|x|<1$, and dividing
by $\pii$, we have
$$V_1=\{(x,y,z,t);\,\max\{|y|,|z|,|t|\}=1,\,|x|<1,\, 
|z^2-t^2+uy^2+u\pii (x/\pii)^2|=1/q\}.$$
The volume of this subset is the volume of $V_0$ of Lemma II.3
multiplied by $1/q$. The lemma follows.
\enddemos
\medskip

\proclaim{II.5. Lemma} When the quadratic form is $x^2-uy^2-uz^2+u\pii t^2$, 
we have
$$\Vol(V_n^0)=\cases
1,              &\qquad\text{if}\ n=0,\\
1/q,            &\qquad\text{if}\ n=1,\\
q^{-n}(1-1/q^2),&\qquad\text{if}\ n\ge 2.
\endcases$$
\endproclaim

\demo{Proof} If $-1\in R^{\times 2}$, the form is
$-u(y^2+z^2-u^{-1}x^2-\pii t^2)$, and its integral has
already been considered in Lemma II.3. Thus we
can take $u=-1$, so the form is
$x^2+y^2+z^2-\pii t^2$. 

In the proof of this lemma we will use Theorem 6.27 of the book 
``Finite Fields'' [LN] by Lidl and Niederreiter. This Theorem 6.27 
asserts that if $f$ is a quadratic form in odd number $n$ of variables 
over the finite field $\F_q$ of $q$ elements, then the number of 
solutions in $\F_q$ of the quadratic equation $f(x_1,x_2,...,x_n)=b$, 
$b\in\F_q$, is $q^{n-1}+q^{(n-1)/2}\eta((-1)^{(n-1)/2}b\det(f))$.
Here $\det(f)$ is the determinant of the symmetric matrix representing 
the quadratic form $f$, and $\eta$ is the quadratic character of $\F_q$:
its value on $\F_q^{\times 2}$ is $1$, on $\F_q^\times-\F_q^{\times 2}$
its value is $-1$, and $\eta(0)=0$.
The case of even $n$ is dealt with in [LN], Theorem 6.26. It asserts 
that -- putting $v(b)=-1$ if $b\not=0$, and $v(0)=q-1$ -- the number of 
solutions of $f=b$ is $q^{n-1}+v(b)q^{(n-2)/2}\eta((-1)^{n/2}\det(f))$; 
but it is not used here. The Theorem 6.27 implies that the equation $f=b$, 
$b=0$, where $f$ is the form $x_1^2+x_2^2+x_3^2$ in $n=3$ variables, has 
$q^2$ solutions over $\F_q$.

{\it Consider the case of} $n=0$. Then
$$V_0=\{(x,y,z,t);\,\max\{|x|,|y|,|z|,|t|\}=1,\, 
|x^2+y^2+z^2-\pii t^2|=1\},$$
namely $V_0=\{(x,y,z,t);\,|x^2+y^2+z^2|=1\}$. 
By Theorem 6.27 of [LN], we have
$$\int\int\int_{|x^2+y^2+z^2|<1}dxdydz=\frac{1}{q}.$$
Hence
$$\Vol(V_0)=\int_{|t|\le 1}\int\int\int_{|x^2+y^2+z^2|=1} dxdydz dt
=1-\int\int\int_{|x^2+y^2+z^2|<1} dxdydz=1-\frac{1}{q}.$$ 

{\it Let us consider the case of} $n\ge 2$. As in Lemma I.0,
recall that any $p$-adic number $a$ such that $|a|\le 1$ 
can be written as a power series in $\pii$:
$$a=\sum_{i=0}^{\infty}a_i\pii^i=a_0+a_1\pii+a_2\pii^2+\dots,\qquad a_i\in R.$$
If $|a|=1/q^n$ we may assume that 
$a_0=a_1=\dots=a_{n-1}=0$ and $a_n\ne 0$. We can write
$$x=\sum_{i=0}^{\infty}x_i\pii^i,\qquad 
y=\sum_{i=0}^{\infty}y_i\pii^i,\qquad 
z=\sum_{i=0}^{\infty}z_i\pii^i,\qquad 
t=\sum_{i=0}^{\infty}t_i\pii^i.$$ 
Their squares are
$$x^2=\sum_{i=0}^{\infty}a_i\pii^i,\qquad 
y^2=\sum_{i=0}^{\infty}b_i\pii^i,\qquad 
z^2=\sum_{i=0}^{\infty}c_i\pii^i,\qquad 
t^2=\sum_{i=0}^{\infty}d_i\pii^i,$$
where
$$a_i=\sum_{j=0}^{i}x_j x_{i-j},\qquad
b_i=\sum_{j=0}^{i}y_j y_{i-j},\qquad
c_i=\sum_{j=0}^{i}z_j z_{i-j},\qquad
d_i=\sum_{j=0}^{i}t_j t_{i-j},$$
and $x_i,\,y_i,\,z_i,\,t_i,\,a_i,\,b_i,\,c_i,\,d_i\in R$.

We have $$x^2+y^2+z^2-\pii t^2=\sum_{i=0}^{\infty}f_i\pii^i,
\qquad f_i\in R,$$
where $f_0=a_0+b_0+c_0$, $f_i=a_i+b_i+c_i-d_{i-1}$ $(i\ge 1)$. 
Since $|x^2+y^2+z^2-\pii t^2|=1/q^n$ we may assume that
$f_0=f_1=...=f_{n-1}=0$ and $f_n\ne 0$. Thus we obtain the relations 
(modulo $\pii$)
$$a_0+b_0+c_0=0,\quad a_i+b_i+c_i-d_{i-1}=0\,(i=1,...,n-1),\quad 
a_n+b_n+c_n-d_{n-1}\ne 0.$$ 

If $a_0=b_0=c_0=0$, it follows that $x_0=y_0=z_0=t_0=0$ (i.e. $|x|<1$, $|y|<1$,
$|z|<1$). Then $a_1=2x_0x_1=0$, $b_1=2y_0y_1=0$, $c_1=2z_0z_1=0$, 
and thus $d_0=a_1+b_1+c_1=0$, i.e. $|t|<1$. This is a
contradiction, since $\max\{|x|,|y|,|z|,|t|\}=1$. 
Assume that $a_0\neq 0$ (i.e. $x_0\neq 0$). From
$a_i+b_i+c_i-d_{i-1}=0$ ($i=1,...,n-1$) it follows that (since $x_0\ne 0$)
$$x_i=(d_{i-1}-b_i-c_i-\sum_{j=1}^{i-1}x_j x_{i-j})/(2x_0),\qquad
x_n\ne (d_{n-1}-b_n-c_n-\sum_{j=1}^{n-1}x_j x_{n-j})/(2x_0),$$
where in the case of $i=1$ the sum over $j$ is empty. 
Thus, we have
$$\int_{|x|=1}\int_{|t|\le 1}\int\int_{|x^2+y^2+z^2-\pi t^2|=q^{-n}}
dydz dt dx=\left(\frac{1}{q}\right)^{n-1}\left(1-\frac{1}{q}\right).$$

By Theorem 6.27 of [LN], we have
$$\int\int\int_{|x^2+y^2+z^2|<1,\,\max\{|x|,|y|,|z|\}=1} dxdydz 
=\frac{1}{q}\left(1-\frac{1}{q^2}\right).$$
Thus 
$$\Vol(V_n)=\frac{1}{q^{n-1}}\left(1-\frac{1}{q}\right)
\times\frac{1}{q^3}\left(q^2-1\right)=\left(1-\frac{1}{q}\right)
\left(1-\frac{1}{q^2}\right)\frac{1}{q^n}.$$

{\it Let us consider the case of} $n=1$. Recall that
$$V_1=\{(x,y,z,t);\,\max\{|x|,|y|,|z|,|t|\}=1,\, 
|x^2+y^2+z^2-\pii t^2|=1/q\}.$$
We consider two cases.
\smallskip
\n (A) Case of $\max\{|x|,|y|,|z|\}<1$, i.e. $|x|<1$,
$|y|<1$, $|z|<1$, and, consequently, $|t|=1$. The contribution
from this case is
$$\int_{|x|<1}\int_{|y|<1}\int_{|z|<1}\int_{|t|=1} dtdzdydx
=\frac{1}{q^3}\left(1-\frac{1}{q}\right).$$

\n (B) Case of $\max\{|x|,|y|,|z|\}=1$. This is the same
as case $n\ge 2$ (with $n=1$). It contributes
$(1-q^{-1})q^{-1}(1-q^{-2})$.

Adding the contributions from Cases (A) and (B) (divided
by $(1-1/q)$), we obtain
$$\Vol(V_1^0)=\frac{1}{q^3}+\frac{1}{q}\left(1-\frac{1}{q^2}\right)
=\frac{1}{q}.$$
The lemma follows.
\enddemos
\medskip

\demo{Proof of Theorem II} We are now ready to complete the
proof of Theorem II. Recall that we need to compute the 
value at $s=0$ ($m=-2$) of $I_s(\r,A,D)$. Here $I_s(\r,A,D)$ coincides 
with the sum $$\sum_{n=0}^{\infty}q^{-nm}\Vol(V_n^0(\r,A,D))$$
which converges for $m>-1$ by Proposition 1 or alternatively
by Lemmas II.1-II.5. The value at $m=-2$ is obtained then by 
analytic continuation of this sum.
\medskip
\n {\it Case of Lemma} II.1. The integral $I_s(\r,A,D)$ is equal to
$$\Vol(V_0^0)+q^{-m}\Vol(V_1^0)+\sum_{n=2}^{\infty}q^{-nm}\Vol(V_n^0)$$
$$=1-\frac{1}{q}+\left(\frac{2}{q}-\frac{1}{q^2}+\frac{1}{q^3}\right) 
\frac{1}{q^m}+2\left(1-\frac{1}{q}\right)q^{-2(m+1)}
\left(1-\frac{1}{q^{m+1}}\right)^{-1}.$$
When $m=-2$, this is 
$$1-\frac{1}{q}+q^2\left(\frac{2}{q}-\frac{1}{q^2}+\frac{1}{q^3}\right)
+2\left(1-\frac{1}{q}\right)\frac{q^2}{1-q}=0.$$

\n {\it Case of Lemma} II.2. The integral $I_s(\r,A,D)$ is equal to
$$\Vol(V_0^0)+q^{-m}\Vol(V_1^0)+\sum_{n=2}^{\infty}q^{-nm}\Vol(V_n^0)$$
$$=1+\frac{1}{q}+\frac{1}{q^2}\left(1-\frac{1}{q}\right) 
\frac{1}{q^m}+\frac{2}{q}\left(1-\frac{1}{q}\right)q^{-2(m+1)}
\left(1-\frac{1}{q^{m+1}}\right)^{-1}.$$
When $m=-2$, this is 
$$1+\frac{1}{q}+1-\frac{1}{q}+\frac{2}{q}\left(1-\frac{1}{q}\right)
\frac{q^2}{1-q}=0.$$

\n {\it Case of Lemmas} II.3, II.4, II.5. The integral $I_s(\r,A,D)$ is 
equal to
$$\Vol(V_0^0)+q^{-m}\Vol(V_1^0)+\sum_{n=2}^{\infty}q^{-nm}\Vol(V_n^0)$$
$$=1+\frac{1}{q}\frac{1}{q^m}+\left(1-\frac{1}{q}\right)
\left(1+\frac{1}{q}\right)q^{-2(m+1)}\left(1-\frac{1}{q^{m+1}}\right)^{-1}.$$
When $m=-2$, this is 
$$1+\frac{1}{q}q^2+\left(1-\frac{1}{q}\right)
\left(1+\frac{1}{q}\right)\frac{q^2}{1-q}=0.$$

The theorem follows.
\enddemos
 
\heading Character computation for type III\endheading

For the $\theta$-conjugacy class of type III we write out the 
representative $g=t\cdot\diag({\pmb r},{\pmb r})$ as
$$
\left(\matrix 
a_1r_1+a_2r_2A & (a_1r_2+a_2r_1)A & (b_1r_1+b_2r_2A)D & (b_1r_2+b_2r_1)AD\\ 
a_1r_2+a_2r_1  & a_1r_1+a_2r_2A   & (b_1r_2+b_2r_1)D  & (b_1r_1+b_2r_2A)D\\
b_1r_1+b_2r_2A & (b_1r_2+b_2r_1)A &  a_1r_1+a_2r_2A   & (a_1r_2+a_2r_1)A\\ 
b_1r_2+b_2r_1  &  b_1r_1+b_2r_2A  &  a_1r_2+a_2r_1    &  a_1r_1+a_2r_2A 
\endmatrix\right).$$
The product ${}^t\v gJ\v$ (where ${}^t\v=(x,y,z,t)$) is equal to
$$(b_1r_2+b_2r_1)(t^2+z^2A-y^2D-x^2AD)
+2(b_1r_1+b_2r_2A)(zt-xyD),$$
where $a_1+a_2\sqrt A\in E_3^\times$ and $b_1+b_2\sqrt A\in E_3^\times$.
The trace is a function of $g$ in the projective group, and 
$r=r_1+r_2\sqrt{A}$ ranges over a set of representatives in $E_3^\times$
($E_3=F(\sqrt{A})$) for $E_3^\times/N_{E/E_3}E^{\times}$.

By definition, the quadratic form can be written as
$${{rb-\tau(rb)}\over{2\sqrt A}}(t^2+z^2A-y^2D-x^2AD)
+(rb+\tau(rb))(zt-xyD).$$
Set $I_s(r,A,D)$ to be equal to
$$\int_{V^0}\left|{{rb-\tau(rb)}\over{2\sqrt A}}(t^2+z^2A-y^2D-
x^2AD)+(rb+\tau(rb))(zt-xyD)\right|^{2(s-1)} dxdydzdt.$$
The property of the numbers $A$, $D$ and $AD$ that we need is that
their square roots generate the three distinct quadratic extensions
of $F$. Thus we may assume that $\{A,D,AD\}=\{u,\pii,u\pii\}$, where 
$u\in R^\times-R^{\times 2}$. Of course with this normalization $AD$
is no longer the product of $A$ and $D$, but its representative in the
set $\{1,u,\pii,u\pii\}$ mod $F^{\times 2}$. 

\proclaim{III.1. Proposition} (i) If $D=u$ and $A=\pii$ (or $\pii u$)
then $\sqrt{A}\not\in N_{E/E_3}E^\times=A^\Z R_3^\times$.\hb
(ii) If $A=u$ and $-1\in R^{\times 2}$, and $D=\pii$ (or $\pii u$) then
$\sqrt{A}\not\in N_{E/E_3}E^\times=(-D)^\Z R_3^{\times 2}$.\hb
(iii) If $A=u=-1\not\in R^{\times 2}$ and $D=\pii$ (or $\pii u$) then
there is $d\in R^\times$ with $d^2+1\in -R^{\times 2}=R^\times-R^{\times 2}$, 
hence $d+i\in R_3^{\times}-R_3^{\times 2}$ ($i=\sqrt{A}$) and so $d+i\in 
E_3^\times-N_{E/E_3}E^\times$.
\endproclaim

\demo{Proof} For (iii) note that $R^\times/\{1+\pii R\}$ is the 
multiplicative group of a finite field $\F$ of $q$ elements. There 
are $1+{1\over 2}(q-1)$ elements in each of the sets $\{1+x^2;\,x\in\F\}$ 
and $\{-y^2;\,y\in\F\}$. As $2(1+{1\over 2}(q-1))>q$, there are $x$, $y$ 
with $1+x^2=-y^2$. But $y\not=0$ as $-1\not\in\F^{\times 2}$. Hence there 
is $x$ with $1+x^2\not\in\F^{\times 2}$, and our $d$ exists.
\enddemos

Since $r$ ranges over a set of representatives for 
$E_3^\times/N_{E/E_3}E^{\times}$, by Proposition III.1
we can choose $br$ to be $1$ or $\sqrt A$ or $d+i$. 
Correspondingly the quadratic form takes one of the three shapes
$$t^2+z^2A-y^2D-x^2AD,\qquad\text{or}\qquad zt-xyD,\qquad 
t^2-z^2-y^2D+x^2D+2d(zt-xyD).$$

Recall that we need to compute
$$\left({\nu\over{\mu}}\right)(\det g){\Delta(g\theta)
\over{\Delta_C(Ng)}}\int_{V^0}|{}^t\v gJ\v|^md\v.\eqno(\ast)$$
Since $\det g=\alpha r\cdot\sigma(\alpha r)\cdot
\tau(\alpha r)\cdot\tau\sigma(\alpha r)$, we have
$$\left({\nu\over{\mu}}\right)(\det g){\Delta(g\theta)\over{\Delta_C(Ng)}}
=|\det g|^{(1-s)/2}\left|{(\alpha r-\sigma(\alpha r))^2\over{\alpha
r\sigma(\alpha r)}}\cdot {\tau(\alpha r-\sigma(\alpha r))^2\over
{\tau(\alpha r)\tau\sigma (\alpha r)}}\right|^{1/2}$$
$$={|4br\tau (br)D|\over |r^2\tau r^2(a^2-b^2D)
(\tau a^2-\tau b^2D)|^{s/2}}.$$
When $s=0$, this is $|br\tau(br)D|$, and $(\ast)$ is independent of $b$.
So we may assume that $b=1$.

\proclaim{III. Theorem} The value of
$|br\tau(br) D|I_s(r,A,D)/(T\phi_0)(\v_0)$ 
at $s=0$ is $2\kappa_{E/E_3}(r)$, where $\kappa_{E/E_3}$ is the nontrivial 
character of $E_3^\times/N_{E/E_3}E^\times$, $E=E_3(\sqrt D)$.
\endproclaim

\demo{Proof} Assume that $br=\sqrt A\not\in N_{E/E_3}E^\times$, thus 
$|br\tau(br)D|=|AD|$, and the quadratic form is $t^2+z^2A-y^2D-x^2AD$. 
If $|A|=1/q$ or $-1$ is a square, we can replace $A$ with $-A$. The 
quadratic form then becomes the same as that of type I. The result of 
the computation is $-2$, see proof of Theorem I, case of anisotropic 
quadratic forms. Since $\kappa_{E/E_3}(\sqrt A)=-1$ we are done in this
case.

If $A=-1$, $br=d+i\not\in N_{E/E_3}E^\times$, the quadratic form is 
$t^2-z^2-y^2D+x^2D+2d(zt-xyD)$. It is equal to $X^2-uY^2-D(Z^2-uT^2)$ 
with $X=t+dz$, $Y=z$, $Z=y+dx$, $T=x$ and $u=d^2+1
\in R^\times-R^{\times 2}$. Since $|D|=1/q$ the quadratic form is 
anisotropic and the result of the computation is $-2$ by the proof of 
Theorem I, case of anisotropic quadratic forms.

Assume that $br=1$, thus $|br\tau(br)D|=|D|$ and the 
quadratic form is $zt-xyD$. Then it is ${1\over 4}$ times
$(z+t)^2-(z-t)^2-D[(x+y)^2-(x-y)^2].$
Since $\max\{|x|,|y|,|z|,|t|\}=1$ implies 
$\max\{|x+y|,|x-y|,|z+t|,|z-t|\}=1$, the result of the computation
is 2 by the proof of Theorem I, cases of Lemmas I.1 and I.3.
The theorem follows.
\enddemos

\heading Character computation for type IV\endheading

For the $\theta$-conjugacy class of type IV we write the representative 
$g=t\cdot\diag({\pmb r},{\pmb r})$ (where $t=h^{-1}t^\ast h$, $t^\ast=
\diag(\alpha,\sigma\alpha,\sigma^3\alpha,\sigma^2\alpha)$) as
$$
\left(\matrix 
a_1r_1+a_2r_2A & (a_1r_2+a_2r_1)A & (b'_1r_1+b'_2r_2A)D &(b'_1r_2+b'_2r_1)AD\\ 
a_1r_2+a_2r_1  & a_1r_1+a_2r_2A   & (b'_1r_2+b'_2r_1)D  &(b'_1r_1+b'_2r_2A)D\\
b_1r_1+b_2r_2A & (b_1r_2+b_2r_1)A &  a_1r_1+a_2r_2A     & (a_1r_2+a_2r_1)A\\ 
b_1r_2+b_2r_1  &  b_1r_1+b_2r_2A  &  a_1r_2+a_2r_1      &  a_1r_1+a_2r_2A 
\endmatrix\right).$$
Here $E_3=F(\sqrt{A})$ is a quadratic extension of $F$ and $E=E_3(\sqrt D)$
is a quadratic extension of $E_3$, thus $A\in F-F^2$ and $D=d_1+d_2\sqrt A\in 
E_3-E_3^2$, $d_i\in F$. 

If $-1\in F^{\times 2}$ we can and do take $D=\sqrt{A}$, where $A$ is a 
nonsquare unit $u$ if $E_3/E$ is unramified, or a uniformizer $\pii$ if 
$E_3/F$ is ramified. If $-1\not\in F^{\times 2}$ and $E_3/F$ is ramified, 
once again we may and do take $A=\pii$ and $D=\sqrt{A}$. 

If $-1\not\in F^{\times 2}$ and $E_3/F$ is unramified, take $A=-1$ and 
note that a primitive 4th root $\zeta=i$ of 1 lies in $E_3$ (and generates 
it over $F$). Then $E/E_3$ is unramified, generated by $\sqrt{D}$, 
$D=d_1+id_2$, and we can (and do) take $d_2=1$ and a unit $d_1=d$ in 
$F^\times$ such that $d^2+1\not\in F^{\times 2}$. Then $D=d+i\not\in 
E_3^{\times 2}$. The existence of $d$ is shown as in the proof of 
Proposition III.1.

Further $\alpha=a+b\sqrt{D}\in E^\times$, where
$a=a_1+a_2\sqrt A\in E_3^\times$, $b=b_1+b_2\sqrt A\in E_3^\times$, and 
$r=r_1+r_2\sqrt{A}\in E_3^\times/N_{E/E_3}E^\times$. The relation
$bD=b'_1+b'_2\sqrt A$ defines $b'_1=b_1d_1+b_2d_2A$ and $b'_2=b_2d_1+b_1d_2$.

The product ${}^t\v gJ\v$ (where ${}^t\v=(x,y,z,t)$) is then equal to
$$(b_1r_2+b_2r_1)(t^2+z^2A)-(b'_1r_2+b'_2r_1)(y^2+x^2A)
+2(b_1r_1+b_2r_2A)zt-2(b'_1r_1+b'_2r_2A)xy.$$

Since $bD=b'_1+b'_2\sqrt A$, this is
$${{br-\sigma(br)}\over{2\sqrt A}}(t^2+z^2A)+(br+\sigma(br))zt$$
$$-{{brD-\sigma(brD)}\over{2\sqrt A}}(y^2+x^2A)-(brD+\sigma(brD))xy.$$

Note that $r$ ranges over a set of representatives for 
$E_3^\times/N_{E/E_3}E^\times$, and $b$ lies in $E_3^\times$.
As $b$ is fixed, we may take $br$ to range over 
$E_3^\times/N_{E/E_3}E^\times$. Thus we may assume that $b=1$.

Further, note that $E_3/F$ is unramified if and only if $E/E_3$ is
unramified. Hence $r$ can be taken to range over $\{1,\pii\}$ if $E_3/F$
is unramified, and over $\{1,u\}$ if $E_3/F$ is ramified,
where $\pii$ is a uniformizer in $F$ and $u$ is a nonsquare
unit in $F$, in these two cases. Thus in both cases we have that
$\sigma(r)=r$, and the quadratic form is equal to the product of $r$ and
$$2zt-{{D-\sigma(D)}\over{2\sqrt A}}(y^2+x^2A)-(D+\sigma(D))xy.$$

Our aim is to compute the value at $s=0$ of the integral $I_s(A,D)$ defined by
$$\int_{V^0}\left|2zt-{{D-\sigma(D)}\over{2\sqrt A}}(y^2+x^2A) 
-(D+\sigma(D))xy\right|^{2(s-1)} dxdydzdt.$$

\proclaim{IV. Theorem} The value of $I_s(A,D)$ at $s=0$ is $0$.
\endproclaim

To prove this theorem we need some lemmas. 

\proclaim{IV.1. Proposition} Up to a change of coordinates, the quadratic form 
$$2zt-{{D-\sigma(D)}\over{2\sqrt A}}(y^2+x^2A)-(D+\sigma(D))xy$$
is equal to either $x^2+\pii y^2-2zt$ or $x^2-uy^2-2zt$ 
with $u\in R^\times-R^{\times 2}$. It is always isotropic.
\endproclaim

\demo{Proof}
In the cases when $D=\sqrt A$, we have $\sigma(D)=-D$. When $D=d+i$,
$\sigma D=d-i$. Thus the quadratic form takes one of the following three shapes
$$2zt-(y^2+\pii x^2),\qquad 2zt-(y^2-u x^2),\qquad 2zt-(y^2-x^2)-2dxy.$$

For the third quadratic form we have
$$2zt-(y^2-x^2)-2dxy\,=\,(x-dy)^2-(d^2+1)y^2+2zt.$$
Recall that $u=d^2+1\in R^\times-R^{\times 2}$. 
After the change of variables $x'=x-dy$, followed by $x'\mapsto x$,
the quadratic form is $x^2-uy^2+2zt$. Change $z\mapsto -z$ to get 
$x^2-uy^2-2zt$.
\enddemos

\proclaim{IV.2. Lemma} When the quadratic form is 
$x^2-uy^2-2zt$, we have
$$\Vol(V_n^0)=\cases
1+1/q^2,                &\qquad\text{if}\ n=0,\\
q^{-n}(1-1/q)(1+1/q^2), &\qquad\text{if}\ n\ge 1.
\endcases$$
\endproclaim

\demo{Proof} {\it Consider the case of} $n=0$. Then
$$V_0=\{(x,y,z,t);\,\max\{|x|,|y|,|z|,|t|\}=1,\, 
|x^2-uy^2-2zt|=1\}.$$

\n (A) Case of $|x^2-uy^2|=1$ and $|zt|<1$. The contribution is the
product of
$$\int\int_{|x^2-uy^2|=1}dxdy=\int_{|x|=1}\int_{|y|\le 1}dydx
+\int_{|x|<1}\int_{|y|=1}dydx=\left(1-\frac{1}{q}\right)
\left(1+\frac{1}{q}\right)$$ 
and
$$\int\int_{|zt|<1}dzdt=\int_{|z|<1}\int_{|t|\le 1}dtdz+
\int_{|z|=1}\int_{|t|<1}dzdt=\frac{1}{q}+\frac{1}{q}
\left(1-\frac{1}{q}\right)
=\frac{1}{q}\left(2-\frac{1}{q}\right).$$

\n (B) Case of $|x^2-uy^2|<1$ and $|zt|=1$ (i.e. $|z|=1$,
$|t|=1$). Since $x^2-u y^2$ does not represent zero non trivially, 
the condition implies that $|x|<1$, $|y|<1$. Thus we obtain
$$\int_{|x|< 1}\int_{|y|< 1}\int_{|z|=1}\int_{|t|=1}
dtdzdydx=\frac{1}{q^2}\left(1-\frac{1}{q}\right)^2.$$

\n (C) Case of $|x^2-uy^2|=1$ and $|zt|=1$. In this
case, once $x$, $y$, and $z$ are chosen, we have that
$|t|=1$, and the condition $|x^2-uy^2-2zt|=1$ implies
$t\not\equiv (x^2-uy^2)/(2z)$ (mod $\pii$). Since
$\int_{|x^2-uy^2|=1}dxdy=1-q^{-2}$, we obtain
$$\int\int_{|x^2-uy^2|=1}\int_{|z|=1}\int_{|t|=1,|x^2-uy^2+zt|=1}
dtdzdxdy=\left(1-\frac{1}{q^2}\right)\left(1-\frac{1}{q}\right)
\left(1-\frac{2}{q}\right).$$

Adding the contributions from Cases (A), (B), and (C) (divided
by $(1-1/q)$), we obtain
$$\Vol(V_0^0)=\frac{1}{q}\left(1+\frac{1}{q}\right) 
\left(2-\frac{1}{q}\right)+\frac{1}{q^2}\left(1-\frac{1}{q}\right) 
+\left(1-\frac{1}{q}\right)\left(1+\frac{1}{q}\right)
\left(1-\frac{2}{q}\right).$$
Once simplified this is equal to $1+1/q^2$.

{\it Consider the case} $n\ge 1$. We have the following two cases.
\smallskip
\n (A) Case of $|z|=1$. Then $x^2-uy^2-2zt=\varepsilon\pii^n$,
where $|\varepsilon|=1$, and $t=(x^2-uy^2-\varepsilon\pii^n)/(2z)$.
Further, $dt=q^{-n}d\varepsilon$, and the contribution from this
case is
$$\int_{|x|\le 1}\int_{|y|\le 1}\int_{|z|=1}\int_{|\varepsilon|=1}
\frac{1}{q^n}d\varepsilon dzdydx=\left(1-\frac{1}{q}\right)
\frac{1}{q^n}\int_{|\varepsilon|=1}d\varepsilon=\frac{1}{q^n}
\left(1-\frac{1}{q}\right)^2.$$  

\n (B) Case of $|z|<1$. If $|t|<1$, then $\max\{|x|,|y|\}=1$,
and since $x^2-uy^2$ does not represent zero non trivially, we
have that $|x^2-uy^2-2zt|=|x^2-uy^2|=1$, which is a contradiction,
since $n\ge 1$. Hence $|t|=1$. We have 
$x^2-uy^2-2zt=\varepsilon\pii^n$, where $|\varepsilon|=1$. Further, from
$$|z|=\left|\frac{x^2-uy^2}{2t}-\frac{\varepsilon}{2t}\pii^n\right|
=|x^2-uy^2-\varepsilon\pii^n|<1,$$ it follows that $|x^2-uy^2|<1$, 
and thus $|x|<1$, $|y|<1$. The contribution from this case is
$$\int_{|x|< 1}\int_{|y|< 1}\int_{|t|=1}\int_{|\varepsilon|=1}
\frac{1}{q^n}d\varepsilon dtdydx=\frac{1}{q^n}
\left(1-\frac{1}{q}\right)^2\frac{1}{q^2}.$$  

Adding the contributions from Cases (A) and (B) (divided
by $(1-1/q)$), we obtain
$$\Vol(V_n^0)=\left(1-\frac{1}{q}\right)\left(1+\frac{1}{q^2}\right) 
\frac{1}{q^n}.$$
The lemma follows.
\enddemos
\medskip

\demo{Proof of Theorem IV} We are now ready to complete the
proof of Theorem IV. Recall that we need to compute the 
value at $s=0$ ($m=-2$) of $I_s(A,D)$. Here $I_s(A,D)$ coincides 
with the sum $$\sum_{n=0}^{\infty}q^{-nm}\Vol(V_n^0(A,D))$$
which converges for $m>-1$. The value at $m=-2$ is obtained then by 
analytic continuation of this sum.
\medskip
\n {\it Case of $x^2+\pii y^2-2zt$}. Make a change of variables 
$z\mapsto 2u^{-1} z'$, followed by $z'\mapsto z$. Thus the quadratic
form is equal to
$$-u^{-1}((z-t)^2-(z+t)^2-u x^2-u\pii y^2).$$
Note that up to a multiple by a unit, this is a form of Lemma II.3. 
Since $\max\{|z|,|t|\}=1$ implies $\max\{|z+t|,|z-t|\}=1$,
the result of that lemma holds for our quadratic form as well.
\medskip
\n {\it Case of $x^2-uy^2-2zt$}. By Lemma IV.2, the integral
$$I_s(A,D)=\Vol(V_0^0)+\sum_{n=1}^{\infty}q^{-nm}\Vol(V_n^0)$$
is equal to
$$1+\frac{1}{q^2}+\left(1-\frac{1}{q}\right)
\left(1+\frac{1}{q^2}\right)q^{-(m+1)}
\left(1-\frac{1}{q^{m+1}}\right)^{-1}.$$
When $m=-2$, this is 
$$1+\frac{1}{q^2}+\left(1-\frac{1}{q}\right)\left(1+\frac{1}{q^2}\right)
\frac{q}{1-q}=0.$$
The theorem follows.
\enddemos
\medskip

\newpage
\def\refe#1#2{\n\hangindent 5em\hangafter1\hbox to 5em{\hfil#1\quad}#2}
\subheading{References}
\medskip

\refe{[BZ]}{I. Bernstein, A. Zelevinsky, Induced representations of 
reductive $p$-adic groups I, {\it Ann. Sci. Ec. Norm. Super.} 10 (1977), 
441-472.}

\refe{[C]}{L. Clozel, Characters of non-connected, reductive $p$-adic groups,
{\it Canad. J. Math.} 39 (1987), 149-167.}

\refe{[DM]}{P. Deligne, J. Milne, Tannakian categories, in {\it Hodge Cycles,
Motives, and Shimura Varieties}, Lecture Notes in Mathematics 900,
Springer-Verlag (1982), 101-228.}

\refe{[Fsym]}{Y. Flicker, On the symmetric-square: Applications of a trace 
formula, {\it Trans. AMS} 330 (1992), 125-152;  Total global comparison, {\it 
J. Funct. Anal.} 122 (1994), 255-278; Unit elements, {\it Pacific J. Math.}
175 (1996), 507-526; {\it Automorphic Representations of Low Rank Groups},
research monograph.}

\refe{[F]}{Y. Flicker, {\it Matching of Orbital Integrals on} GL(4) {\it and} 
GSp(2), Memoirs AMS 137 (1999), 1-114.}

\refe{[F$'$]}{Y. Flicker, {\it Lifting Automorphic Forms of} PGSp(2) {\it and} 
SO(4) {\it to} PGL(4), research monograph.}

\refe{[FK]}{Y. Flicker, D. Kazhdan, On the symmetric-square. Unstable
local transfer, {\it Invent. Math.} 91 (1988), 493-504.}

\refe{[FZ]}{Y. Flicker, D. Zinoviev, On the symmetric-square. Unstable
twisted characters, {\it Israel J. Math.} 134 (2003), 307-316.}

\refe{[FZ$'$]}{Y. Flicker, D. Zinoviev, Twisted character of a small
representation of GL(4), preprint.}
 
\refe{[H]}{Harish-Chandra, Admissible invariant distributions on reductive 
$p$-adic groups, {\it Queen's papers in Pure and Applied Math.} 48 (1978), 
281-347.}

\refe{[K]}{D. Kazhdan, On liftings, in {\it Lie Groups Representations II},
Springer Lecture Notes on Mathematics 1041 (1984), 209-249.}

\refe{[KS]}{R. Kottwitz, D. Shelstad, {\it Foundations of Twisted Endoscopy}, 
Asterisque 255 (1999), vi+190 pp.}
 
\refe{[LN]}{R. Lidl, H. Niederreiter, {\it Finite Fields}, Cambridge Univ.
Press, 1997.}

\refe{[Wa]}{J.-L. Waldspurger, Sur les int\'egrales orbitales tordues pour 
les groupes lin\'eaires: un lemme fondamental, {\it Canad. J. Math.} 43 
(1991), 852-896.}

\refe{[W]}{R. Weissauer, A special case of the fundamental lemma, preprint.}


\enddocument